\documentclass[a4paper,11pt]{article}
\usepackage{authblk}
\usepackage[british]{babel}
\usepackage[mono=false]{libertine}
\usepackage{setspace}
\usepackage[utf8]{inputenc} 
\usepackage[T1]{fontenc}    
\usepackage{url}            
\usepackage{microtype}
\usepackage{graphicx}
\usepackage{amssymb}
\usepackage{enumitem}
\usepackage{subcaption}
\usepackage{booktabs} 
\usepackage{amsfonts}       
\usepackage{amsmath,amssymb,amsthm,mathrsfs,bbm,bm}
\usepackage{physics}
\usepackage{color}
\usepackage[usenames, dvipsnames]{xcolor}
\usepackage{enumitem}

\usepackage{algorithm}
\usepackage{algorithmic}

\usepackage[top=2.5cm, bottom=2.5cm, left=2cm, right=2cm]{geometry}
\usepackage[pdfpagemode=UseNone,bookmarks=false, bookmarksopen=false,colorlinks=true,urlcolor=RoyalBlue,citecolor=YellowOrange,citebordercolor=RoyalBlue,linkcolor=RoyalBlue]{hyperref}

\def\new{\text{new}}

\newcommand{\extr}{\textrm{\textbf{extr}}}
\newcommand{\argmin}{\textrm{argmin}}

\def\dd{\text{d}}

\def\sign{\text{sign}}
\def\erf{\text{erf}}

\def\mat#1{\text{#1}}
\renewcommand{\vec}[1]{\bm{#1}}

\newtheorem{theorem}{Theorem}

\numberwithin{equation}{section}

\begin{document}

\title{Generalisation error in learning with random features\\ 
and the hidden manifold model}
\author{Federica Gerace$^\dagger$, Bruno Loureiro$^\dagger$, \\ 
Florent Krzakala$^\star$, Marc M{\'e}zard $^{\star}$, Lenka Zdeborov{\'a}$^\dagger$}
\date{
$\dagger$ \textit{Institut de Physique Th\'eorique \\
CNRS \& CEA \& Universit\'e Paris-Saclay, Saclay, France}\\
$\star$ \textit{LPENS, CNRS \& Sorbonnes Universit´es,  \\ 
Ecole Normale Sup´erieure, PSL University, Paris, France}\\
\vspace{1cm}
}

\maketitle

\begin{abstract}
We study generalised linear regression and classification for a synthetically generated dataset encompassing different problems of interest, such as learning
with random features, neural networks in the lazy training regime, and the hidden manifold model.  We consider the high-dimensional regime and using the replica
method from statistical physics, we provide a closed-form expression for the asymptotic generalisation performance in these problems, valid in both the under- and over-parametrised regimes and for a broad choice of generalised linear model loss functions. In particular, we show how to obtain analytically the so-called double descent behaviour for logistic regression with a peak at the interpolation threshold, we illustrate the superiority of orthogonal against random Gaussian projections in learning with random features, and discuss the role played by correlations in the data generated by the hidden manifold model. Beyond the interest in these particular problems, the theoretical formalism introduced in this manuscript provides a path to further extensions to more complex tasks.
\end{abstract}

\newpage
\makeatletter
\def\l@subsubsection#1#2{}
\makeatother

{
  \hypersetup{linkcolor=black}
\begin{spacing}{0.9}
  \tableofcontents
  \end{spacing}
}

\newpage
\section{Introduction}
One of the most important goals of learning theory is to provide
generalisation bounds describing
the quality of learning a given task
as a function of the number of samples. Existing results fall short of
being directly relevant for the state-of-the-art deep learning methods
\cite{zhang2016understanding,neyshabur2017exploring}. Consequently,
providing tighter
results on the generalisation
error is currently a very active research subject. 
The traditional learning theory approach to
generalisation follows for instance the
Vapnik-Chervonenkis \cite{vapnik1998statistical} or Rademacher
\cite{bartlett2002rademacher} worst-case type bounds, and many of
their more recent extensions \cite{shalev2014understanding}. An alternative approach,
followed also in this paper, has been
pursued for decades, notably in statistical physics, where the generalisation
ability of neural networks was analysed for a range of 
``typical-case'' scenario {\it for synthetic data arising from a
  probabilistic model}
\cite{seung1992statistical,watkin1993statistical,
advani2013statistical,advani2017high,aubin2018committee,
candes2020, hastie2019surprises,mei2019generalisation,goldt2019modelling}. While
at this point it is not clear which approach
will lead to a complete
generalisation theory of deep learning, it is worth
pursuing both directions. 

The majority of works following the statistical physics approach study the generalisation error in the
so-called teacher-student framework, where the input data are element-wise
i.i.d. vectors, and the labels are generated by a teacher neural
network. In contrast, in most of real scenarios the input data do not span uniformly the input space, but rather live close to a lower-dimensional
manifold.
The traditional focus onto i.i.d. Gaussian input vectors is an important
limitation that has been recently stressed in
\cite{MezardHopfield,goldt2019modelling}. In
\cite{goldt2019modelling}, the authors proposed a model
of synthetic data to mimic the latent structure of real data, named the 	\emph{hidden manifold model}, and analysed the
learning curve of one-pass stochastic gradient descent algorithm in a two-layer neural network with
a small number of hidden units also known as committee machine.

Another key limitation of the majority of existing works stemming from statistical physics is
that the learning curves were only computed for neural
networks with a few hidden units. In particular, the input dimension
is considered large, the number of samples is a constant times the
input dimension and the number of hidden units is of order one. Tight
learning curves were only very recently analysed for two-layer neural
networks with more hidden units. These studies addressed in particular
the case of networks that have a fixed first layer with random i.i.d. Gaussian weights
\cite{hastie2019surprises,mei2019generalisation}, or the lazy-training regime where the individual weights
change only infinitesimally during training, thus not learning any
specific features \cite{chizat2018note,jacot2018neural,geiger2019disentangling}. 

In this paper we compute the generalisation error and the corresponding
learning curves, i.e. the test error as a function of the number of
samples for a model of high-dimensional data that encompasses at least the following cases:
\begin{itemize}
      \item{generalised linear regression and classification for data
        generated by the hidden manifold model (HMM) of
        \cite{goldt2019modelling}. The HMM can be seen as a
        single-layer
        generative neural network with i.i.d. inputs and a rather
        generic feature matrix \cite{louart2018random,
  goldt2019modelling}. }  
\item{Learning data generated by the teacher-student model with a
    random-features neural network \cite{NIPS2007_3182}, with a very
    generic feature matrix, including deterministic ones. This
    model is also interesting because of its connection with the lazy
    regime, that is equivalent to the random features model with
    slightly more complicated
    features~\cite{chizat2018note,hastie2019surprises,mei2019generalisation}.}
\end{itemize}

We give a closed-form expression for the generalisation error in the
high-dimensional limit, obtained using a non-rigorous heuristic method
from statistical physics known as the replica method
\cite{mezard1987spin}, that has already shown its remarkable efficacy
in many problems of machine
learning~\cite{seung1992statistical,engel2001statistical,advani2013statistical,zdeborova2016statistical}, with many of its predictions being rigorously proven, e.g.
\cite{talagrand2006parisi,barbier2019optimal}. While
in the present model it remains an open problem to derive a rigorous
proof for our results, we shall provide numerical support that the
formula is indeed exact in the high-dimensional limit, and extremely
accurate even for moderately small system sizes.

\subsection{The model}
\label{sec:model}
%
We study high-dimensional regression and classification
for a {\it synthetic} dataset
${\cal D} = \{(\vec{x}^{\mu}, y^{\mu})\}_{\mu=1}^{n}$ where each sample $\mu$ is
created in the following three steps: (i) First, for each sample $\mu$ we create a vector
$\vec{c}^{\mu}\in\mathbb{R}^{d}$ as
\begin{equation}
  \vec{c}^{\mu} \sim \mathcal{N}(0,\mat{I}_{d})\, ,
  \label{model-def-c}
\end{equation}
(ii) We then draw $\vec{\theta}^{0}\in{\mathbb R}^{d}$ 
from a separable distribution $P_{\theta}$ and draw independent labels $\{y^\mu\}_{\mu=1}^n$ from a (possibly probabilistic) rule $f^{0}$: 
\begin{align}
  y^{\mu} = f^{0}\left(\frac{1}{\sqrt{d}}\vec{c}^{\mu}\cdot
  \vec{\theta}^{0}\right)\in\mathbb{R} \, .
  \label{model-def-y}
\end{align}
(iii) The input data
points $\vec{x}^{\mu}\in\mathbb{R}^{p}$ are created by a one-layer
generative network with fixed and normalised weights $\mat{F} \in {\mathbb R}^{d
  \times p}$ and an activation function $\sigma:\mathbb{R}\to\mathbb{R}$, acting component-wise:
\begin{align}
\vec{x}^{\mu} = \sigma\left(\frac{1}{\sqrt{d}}  \mat{F}^{\top}  \vec{c}^{\mu}
  \right) \, .
\label{model-def-x}
\end{align}
We study the problem of supervised learning for the dataset~${\cal D}$
aiming at achieving a low
generalisation error $\epsilon_{\rm g}$ on a new sample
${\vec{x}^{\rm new},y^{\rm new}}$ drawn by the same rule as above, where:
\begin{align}
\epsilon_{\rm g} = \frac{1}{4^k}\mathbb{E}_{{\vec{x}^{\rm new},y^{\rm new}}}
  \left[\left(\hat{y}_{\vec{w}}(\vec{x}^{\rm new}) - y^{\new}
  \right)^2\right].
  \label{def-gen}
\end{align}
with $k = 0$ for regression task and $k=1$ for classification task. Here, $\hat{y}_{\vec{w}}$ is the prediction on the new label
$y^{\new}$ of the form:
\begin{align}
  \hat{y}_{\vec{w}}(\vec{x}) =\hat  f\left(\vec{x}\cdot \hat{\vec{w}}\right).
  \label{def-f-fit}
\end{align}
The weights $\hat{\vec{w}}\in\mathbb{R}^{p}$ are learned by
minimising a loss function with a ridge regularisation term (for
$\lambda \ge 0$) and defined as 
\begin{align}
	\hat{\vec{w}} = \underset{\vec{w}}{\argmin}~\left[ 
  \sum\limits_{\mu=1}^{n}\ell(y^{\mu}, \vec{x}^\mu\cdot
  \vec{w})+\frac{\lambda}{2}||\vec{w}||_{2}^2 \right],
    \label{eq-loss}
\end{align}
\noindent
where $\ell(\cdot,\cdot)$ can be, for instance, a logistic, hinge, or square
loss. Note that although our formula is valid for any $f^{0}$ and $\hat{f}$, we take $f^{0}=\hat{f}=\sign$, for the classification tasks and $f^{0} =\hat{f} = \text{id}$ for the regression tasks studied here. We now describe in more detail the above-discussed reasons why
this model is of interest for machine learning.
\paragraph{Hidden manifold model:} The dataset
${\cal D}$ can be seen as
generated from the \emph{hidden manifold model} introduced in 
\cite{goldt2019modelling}. From this perspective, although
$\vec{x}^{\mu}$ lives in a $p$ dimensional space, it is parametrised
by a latent $d$-dimensional subspace spanned by the rows of the matrix~$\mat{F}$
which are "hidden" by the application of a scalar non-linear function $\sigma$. The labels ${y}^{\mu}$ are drawn from a
generalised linear rule defined on the latent $d$-dimensional
subspace via eq.~(\ref{model-def-y}). In modern machine learning parlance, this can be seen as data generated by a
one-layer generative neural network, such as those trained by generative
adversarial networks or variational auto-encoders with random 
Gaussian inputs ${\bf c}^\mu$ and a rather generic weight matrix $\mat{F}$
\cite{GAN,VAE,louart2018random,el2020random}. 

\paragraph{Random features:} The model considered in this paper is also an
instance of the random features learning discussed in
\cite{NIPS2007_3182} as a way to speed up
kernel-ridge-regression. From this perspective,
the $\vec{c}^{\mu}$s $\in\mathbb{R}^{d}$ are regarded as a set of
$d$-dimensional i.i.d. Gaussian data points, which are projected by a
feature matrix
$\mat{F} = (\vec{f}_{\rho})_{\rho=1}^{p}\in\mathbb{R}^{d\times p}$
into a higher dimensional space, followed by a
non-linearity~$\sigma$. In the $p\to\infty$ limit of infinite number of features, performing regression on ${\cal D}$ is equivalent to
kernel regression on the ${\bf c}^{\mu}$s with a deterministic kernel
$K({\bf c}^{\mu_1}, {\bf c}^{\mu_2}) = {\mathbb{E}}_{\bf f}\left[
  \sigma({\bf f}\cdot {\bf c}^{\mu_1}/\sqrt{d} ) \cdot \sigma({\bf f}\cdot
  {\bf c}^{\mu_2}/\sqrt{d} )\right]$ where ${\bf f} \in {\mathbb R}^d$ is
sampled in the same way as the rows of $\mat{F}$. Random
features are also intimately linked with the lazy training regime,
where the weights of a neural network stay close to their initial
value during training. The training is lazy as opposed to a ``rich''
one where the weights change enough to learn useful features. In this
regime, neural networks become equivalent to a random feature model
with correlated features
\cite{chizat2018note,du2018gradient,allen2018convergence,woodworth2019kernel,jacot2018neural,geiger2019disentangling}.

\subsection{Contributions and related work}
%
The main contribution of this work is a closed-form expression for the
generalisation error $\epsilon_{\rm g}$, eq.~(\ref{eq:eg_main}), that is
valid in the high-dimensional limit where the number of samples~$n$,
and the two dimensions $p$ and $d$ are large, but their respective
ratios are of order one, and for generic sequence of matrices
$\mat{F}$ satisfying the following \emph{balance conditions}:
\begin{align}
\frac{1}{\sqrt{p}}\sum\limits_{i=1}^{p}w_{i}^{a_{1}}w_{i}^{a_{2}}\cdots w_{i}^{a_{s}} \mat{F}_{i\rho_{1}}\mat{F}_{i\rho_{2}}\cdots \mat{F}_{i\rho_{q}} = O(1), 
\label{eq:balance}
\end{align}
\noindent \noindent where $\{\vec{w}^{a}\}_{a=1}^{r}$ are $r$
independent samples from
the Gibbs measure (\ref{eq:main:Gibbs_definition}),
and
$\rho_{1}, \rho_{2}, \cdots, \rho_{q}\in\{1,\cdots, d\}$,
$a_{1}, a_{2}, \cdots, a_{s}\in\{1,\cdots, r\}$ are an arbitrary
choice of subset of indices, with $s, q\in\mathbb{Z}_{+}$. 
The non-linearities
$f^{0}, \hat f, \sigma$ and the loss function $\ell$ can be
arbitrary.
Our result for the generalisation error stems from the replica method and we
conjecture it to be exact for convex
loss functions $\ell$. It can also be
useful for non-convex loss functions but in those cases it is possible
that the so-called replica symmetry breaking \cite{mezard1987spin}
needs to be taken into account to obtain an exact expression. In the
present paper we hence focus on convex loss functions $\ell$ and leave
the more general case for future work.
The final formulas are simpler for nonlinearities $\sigma$ that give
zero when integrated over a centred Gaussian variable, and we hence
focus on those cases. 

An interesting application of our setting is ridge regression, i.e. 
taking $\hat{f}(x)\!=\!x$ with square loss, and random i.i.d. Gaussian feature matrices. For this
particular case \cite{mei2019generalisation} proved an equivalent
expression. Indeed, in this case there is an explicit solution of eq.~\eqref{eq-loss} that can be rigorously studied with random matrix theory. In a subsequent work \cite{montanari2019generalisation}
derived heuristically a formula for the special case of random i.i.d. Gaussian feature matrices for the maximum margin classification, corresponding to the
hinge loss function in our setting, with the difference,
however, that the labels $y^{\mu}$ are generated from the ${\bf x}^{\mu}$
instead of the variable $\bf c^{\mu}$ as in our case.

Our main technical contribution is thus to provide a generic formula
for the model described in Section \ref{sec:model} for any loss function and for
fairly generic features $\mat{F}$, including for instance deterministic ones.

The authors of \cite{goldt2019modelling} analysed the learning dynamics of a neural network containing
several hidden units using a one-pass stochastic gradient descent
(SGD) for exactly the same model of data as here. In this online setting, the algorithm is never exposed to a sample twice, greatly simplifying the analysis as what has been learned at a given epoch can be considered independent of the randomness of a new sample. Another motivation of the present work is thus to study the sample complexity for this model (in our case only a bounded number of samples is available, and the one-pass SGD would be highly suboptimal).

An additional technical contribution of our work is to derive an extension of the equivalence between the considered data model and a model with Gaussian
covariate, that has been observed and conjectured to hold rather
generically in both \cite{goldt2019modelling,montanari2019generalisation}. While we do not provide a rigorous proof for this equivalence, we show that it arises naturally using the replica method, giving further evidence for its validity.

Finally, the analysis of our formula for particular machine learning tasks of interest allows for an analytical investigation of a rich phenomenology that is also observed empirically in real-life scenarios. In particular
\begin{itemize}
\item The double descent behaviour, as termed in
  \cite{belkin2018reconciling} and exemplified in
  \cite{spigler2018jamming}, is exhibited for the non-regularized
  logistic regression loss. The peak of worst generalisation does not
  corresponds to $p=n$ as for the square loss
  \cite{mei2019generalisation}, but rather corresponds to the
  threshold of linear separability of the dataset.
   We also characterise the location of this threshold,
  generalising the results of \cite{candes2020} to our model.
\item When using projections to approximate kernels, it has been observed
  that orthogonal features ${\mat{F}}$ perform better than random i.i.d.
  \cite{choromanski2017unreasonable}. We show that this behaviour
  arises from our analytical formula, illustrating the "unreasonable
  effectiveness of structured random orthogonal embeddings''\cite{choromanski2017unreasonable}.
\item We compute the phase diagram for the generalisation error for the
  hidden manifold model and discuss the dependence on the various
  parameters, in particular the ratio between the ambient and latent
  dimensions.
\end{itemize}
\section{Main analytical results}
\label{main-sec:main:replica_result}
We now state our two main analytical results.
The replica computation used here is in spirit
similar to the one performed in a number of tasks for linear and
generalised linear
models \cite{gardner1989three,seung1992statistical,kabashima2009typical,krzakala2012statistical},
but requires a significant extension to account for the structure of
the data.  We refer the reader to the supplementary material
Sec.~\ref{sec:app:replicas} for the detailed and lengthy
derivation of the final formula. The resulting expression is
conjectured to be exact and, as we shall see, observed to be accurate
even for relatively small dimensions in simulations. Additionally,
these formulas reproduce the rigorous results of
\cite{mei2019generalisation}, in the simplest particular case of a Gaussian
projection matrix and ridge regression task. It remains a challenge to
prove them rigorously in broader generality.
\begin{figure}[t]
  \begin{center}
  \centerline{\includegraphics[width=0.5\columnwidth]{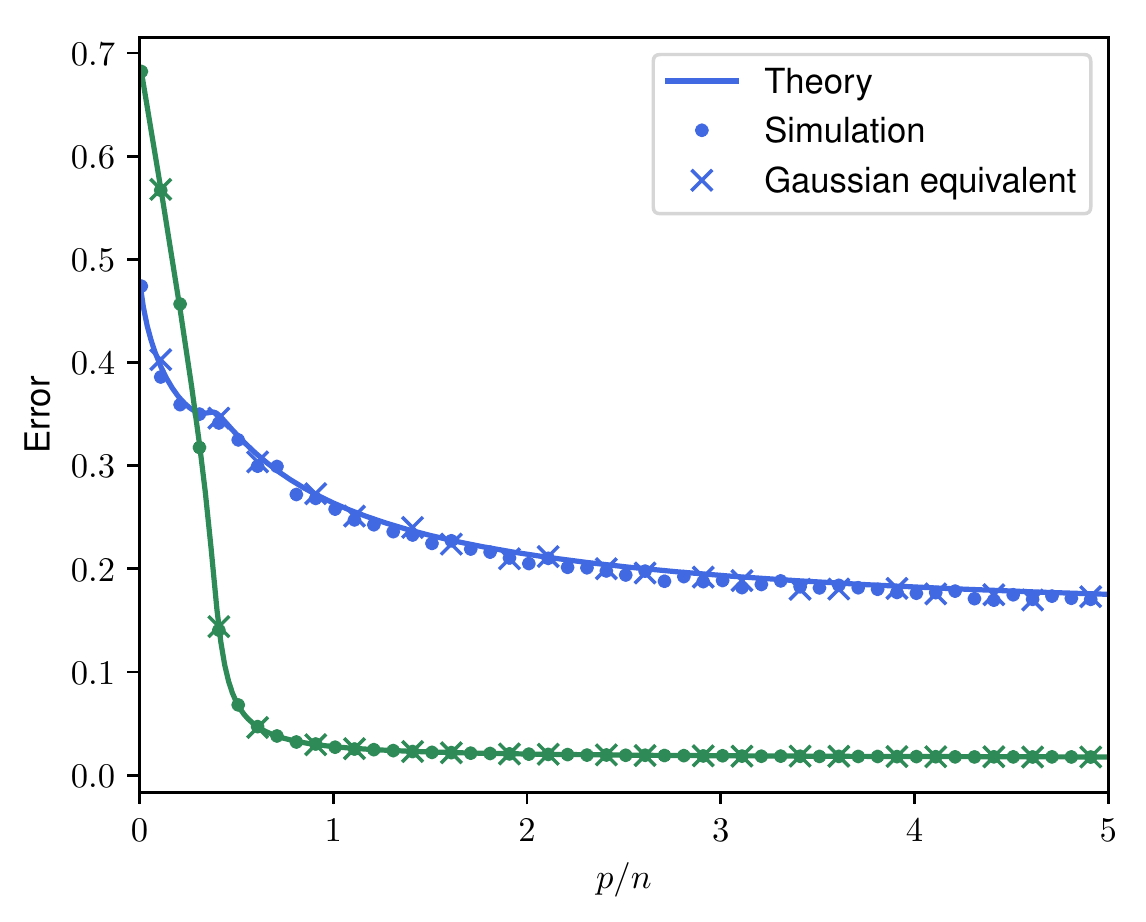}}
  \vspace{-0.6cm}
  \caption{Comparison between theory (full line), and simulations with
    dimension $d=200$ on the original model (dots),
    eq.~(\ref{model-def-x}), with $\sigma = \sign$, and the Gaussian
    equivalent model (crosses), eq.~(\ref{model-def-x-GET}), for
    logistic loss, regularisation $\lambda \!=\! 10^{-3}$, $n/d =
    3$. Labels are generated as
    $y^{\mu}\!=\!\sign\left(\vec{c}^{\mu}\cdot\vec{\theta}^{0}\right)$ and $\hat{f}=\sign$. Both
    the training loss (green) and generalisation error (blue) are
    depicted. The theory and the equivalence with the Gaussian model
    are observed to be very accurate even at dimensions as small as
    $d=200$.}
\label{fig_comparison}
\end{center}
\end{figure}

\subsection{Generalisation error from replica method}
\label{sec:main:replica_result}
%
Let $\mat{F}$ be a feature matrix satisfying the balance condition stated in eq.~\eqref{eq:balance}. Then, in the high-dimensional limit where $p,d,n \to \infty$ with $\alpha =
n/p$, $\gamma = d/p$ fixed, the generalisation error,
eq.~(\ref{def-gen}), of the model introduced in Sec.~(\ref{def-gen})
for $\sigma$ such that its integral over a centered Gaussian variable
is zero (so that $\kappa_0=0$ in eq.~(\ref{model-def-x-GET})) is given by the following easy-to-evaluate integral:
\begin{align}
\lim_{n \to \infty}\epsilon_{g} = \mathbb{E}_{\lambda,\nu}
  \left[(f^{0}(\nu)-\hat f(\lambda))^2\right] \, ,\label{eq:eg_main}
\end{align}
\noindent where $f^{0}(.)$ is defined in (\ref{model-def-y}), $\hat
f(.)$ in  (\ref{def-f-fit})
and $(\nu,\lambda)$ are jointly Gaussian random variables with zero mean
and covariance matrix:
\begin{align}
	\Sigma = 
	\begin{pmatrix}
		\rho &M^{\star} \\
		M^{\star} & Q^{\star}\end{pmatrix}\in\mathbb{R}^2
\end{align}
\noindent with $M^{\star} =
\kappa_{1}m_{s}^{\star}$, $Q^{\star} =
\kappa_{1}^2q^{\star}_{s}+\kappa_{\star}^2q_{w}^{\star}$. The
constants $\kappa_{\star},\kappa_{1}$ depend on the nonlinearity
$\sigma$ via eq.~(\ref{model-def-x-GET}), and $q^{\star}_{s}, q_{w}^{\star}, m_{s}^{\star}$, defined as:
\begin{align}
\rho = \frac{1}{d}||\vec{\theta}^{0}||^2 && q^{\star}_{s} = \frac{1}{d}\mathbb{E}||\mat{F}\vec{\hat{w}}||^2 && q^{\star}_{w} = \frac{1}{p} \mathbb{E}||\vec{\hat{w}}||^2 && m^{\star}_{s} = \frac{1}{d}\mathbb{E}\left[(\mat{F}\hat{\vec{w}})\cdot \vec{\theta}^{0}\right]
\label{eq:optimal_overlaps}
\end{align}
The values of these parameters correspond to the solution of the
optimisation problem in eq.~\eqref{eq-loss}, and can be obtained as the fixed point solutions of the following set of self-consistent \emph{saddle-point equations}:
\begin{align}
	&\begin{cases}
		\hat{V}_{s} = \frac{\alpha\kappa_{1}^{2}}{\gamma V} \mathbb{E}_{\xi}\left[\int_{\mathbb{R}}\dd y~\mathcal{Z}\left(y,\omega_{0}\right)\left(1 - \partial_{\omega}\eta\left(y,\omega_{1}\right) \right)\right],\\
		\hat{q}_{s} = \frac{\alpha\kappa_{1}^2}{\gamma V^2} \mathbb{E}_{\xi}\left[\int_{\mathbb{R}}\dd y~\mathcal{Z}\left(y,\omega_{0}\right)\left(\eta\left(y,\omega_{1}\right)-\omega_{1}\right)^2\right],\\
		\hat{m}_{s} = \frac{\alpha\kappa_{1}}{\gamma V} \mathbb{E}_{\xi}\left[\int_{\mathbb{R}}\dd y~\partial_{\omega}\mathcal{Z}\left(y,\omega_{0}\right)  (\eta\left(y,\omega_{1}\right)-\omega_{1})\right],\\ \\
		\hat{V}_{w} = \frac{\alpha\kappa_{\star}^2}{V} \mathbb{E}_{\xi}\left[\int_{\mathbb{R}}\dd y~\mathcal{Z}\left(y,\omega_{0}\right)\left(1 -\partial_{\omega}\eta\left(y,\omega_{1}\right)\right)\right],\\
		\hat{q}_{w} = \frac{\alpha \kappa_{\star}^2}{V^2}\mathbb{E}_{\xi}\left[\int_{\mathbb{R}}\dd y~\mathcal{Z}\left(y,\omega_{0}\right) \left(\eta\left(y,\omega_{1}\right)-\omega_{1}\right)^2\right],
	\end{cases}&&
	\begin{cases}
	V_{s} = \frac{1}{\hat{V}_{s}}\left(1-z~g_{\mu}(-z)\right), \\
	q_{s} = \frac{\hat{m}_{s}^2+\hat{q}_{s}}{\hat{V}_{s}^2}\left[1-2z g_{\mu}(-z)+z^2 g_{\mu}'(-z)\right]\\ \qquad-\frac{\hat{q}_{w}}{(\lambda+\hat{V}_{w})\hat{V}_{s}}\left[-z g_{\mu}(-z)+z^2 g'_{\mu}(-z)\right], \\
	m_{s} = \frac{\hat{m}_{s}}{\hat{V}_{s}}\left(1-z~g_{\mu}(-z)\right),	\\ \\
	V_{w} =\frac{\gamma}{\lambda+\hat{V}_{w}}\left[\frac{1}{\gamma}-1+z g_{\mu}(-z)\right], \\
	q_{w} = \gamma\frac{\hat{q}_{w}}{(\lambda+\hat{V}_{w})^2}\left[\frac{1}{\gamma}-1+z^2 g'_{\mu}(-z)\right],\\
	\qquad-\gamma\frac{\hat{m}_{s}^2+\hat{q}_{s}}{(\lambda+\hat{V}_{w})\hat{V}_{s}}\left[-z g_{\mu}(-z)+z^2g'_{\mu}(-z)\right],\label{eq:saddle_point}
\end{cases}
\end{align}
\noindent written in terms of the following auxiliary variables
$\xi\sim\mathcal{N}(0,1)$, $z = \frac{\lambda+\hat{V}_{w}}{\hat{V}_{s}}$ and functions:
\begin{align}
\eta(y, \omega) &= \underset{x\in\mathbb{R}}{\argmin}\left[\frac{(x-\omega)^2}{2 V}+\ell(y,x)\right],\notag	\\
\mathcal{Z}(y, \omega) &= \int\frac{\dd x}{\sqrt{2\pi V^{0}}}e^{-\frac{1}{2V^{0}}\left(x-\omega\right)^2}\delta\left(y-f^{0}(x)\right)
\label{eq:etaZ}
\end{align}

\noindent where $V= \kappa_{1}^{2}V_{s}+\kappa_{\star}^2 V_{w}$,
$V^{0} = \rho-\frac{M^2}{Q}$, $Q =
\kappa_{1}^2q_{s}+\kappa_{\star}^2q_{w}$, $M = \kappa_{1}m_{s}$, $\omega_{0}=\left(M/\sqrt{Q}\right)\xi$ and $\omega_{1}=\sqrt{Q}\xi$. In
the above, we assume that the matrix
$\mat{F}\mat{F}^{\top}\in\mathbb{R}^{d\times d}$ associated to the
feature map $\mat{F}$ has a well behaved spectral density, and denote
$g_{\mu}$ its Stieltjes transform.
 
The training loss on the dataset ${\cal D} =
\{\vec{x}^{\mu},y^{\mu}\}_{\mu=1}^{n}$ can also be obtained from the
solution of the above equations as
 \begin{align}
 \lim\limits_{n\to\infty}\epsilon_{t}= \frac{\lambda}{2\alpha}q_{w}^{\star}+\mathbb{E}_{\xi, y}\left[\mathcal{Z}\left(y,\omega_{0}^{\star}\right)\ell\left(y,\eta(y,\omega_{1}^{\star})\right)\right]
 \end{align}
\noindent where as before $\xi\sim\mathcal{N}(0,1)$, $y\sim
\text{Uni}(\mathbb{R})$ and $\mathcal{Z},\eta$ are the same as in
eq.~\eqref{eq:etaZ}, evaluated at the solution of the above
saddle-point equations $\omega_{0}^{\star} =
\left(M^{\star}/\sqrt{Q^{\star}}\right)\xi$,
$\omega_{1}^{\star}=\sqrt{Q^{\star}}\xi$.

%
\paragraph{Sketch of derivation ---}
We now sketch the derivation of the above result. A complete and detailed account can be found in Sec.~\ref{sec:app:replicas} of the supplementary material. 
The derivation is based on the key observation that in the high-dimensional limit the asymptotic generalisation error only depends on the solution $\hat{\vec{w}}\in\mathbb{R}^{p}$ of eq.~\eqref{def-f-fit} through the scalar parameters $(q_{s}^{\star}, q_{w}^{\star}, m_{s}^{\star})$ defined in eq.~\eqref{eq:optimal_overlaps}. The idea is therefore to rewrite the high-dimensional optimisation problem in terms of only these scalar parameters.

The first step is to note that the solution of eq.~\eqref{eq-loss} can be written as the average of the following Gibbs measure
\begin{align}
 \pi_{\beta}(\vec{w}|\{\vec{x}^{\mu},y^{\mu}\})	 = \frac{1}{\mathcal{Z}_{\beta}}e^{-\beta\left[\sum\limits_{\mu=1}^{n}\ell\left(y^{\mu}, \vec{x}^{\mu}\cdot \vec{w}\right) + \frac{\lambda}{2}||\vec{w}||_{2}^2\right]},
  \label{eq:main:Gibbs_definition}
\end{align}
\noindent in the limit $\beta\to\infty$. Of course, we have not gained much, since an exact calculation of $\pi_{\beta}$ is intractable for large values of $n, p$ and $d$. This is where the replica method comes in. It states that the distribution of the free energy density $f = -\log\mathcal{Z}_{\beta}$ (when seen as a random variable over different realisations of dataset $\mathcal{D}$) associated with the measure $\mu_{\beta}$ concentrates, in the high-dimensional limit, around a value $f_{\beta}$ that depends only on the averaged \emph{replicated partition function} $\mathcal{Z}_{\beta}^{r}$ obtained by taking $r>0$ copies of $\mathcal{Z}_{\beta}$:
\begin{align}
f_\beta= \lim\limits_{r\to 0^{+}}\frac{\dd}{\dd r}\lim\limits_{p\to\infty}\left[ -\frac{1}{p} \left(\mathbb{E}_{\lbrace\vec{x}^{\mu},y^{\mu}\rbrace}\mathcal{Z}_{\beta}^{r}\right)\right].\label{eq:free_en}
 \end{align} 
Interestingly, $\mathbb{E}_{\lbrace\vec{x}^{\mu},y^{\mu}\rbrace}\mathcal{Z}_{\beta}^{r}$ can be
computed explicitly for $r\in\mathbb{N}$, and the limit $r\to 0^{+}$ is taken by analytically continuing to $r>0$ (see Sec.~\ref{sec:app:replicas} of the
supplementary material). The upshot is that $\mathcal{Z}^{r}$ can be written as
\begin{align}
	\mathbb{E}_{\lbrace\vec{x}^{\mu},y^{\mu}\rbrace}\mathcal{Z}_{\beta}^{r} \propto \int\dd q_{s}\dd q_{w}\dd m_{s} ~e^{p\Phi_{\beta}^{(r)}(m_{s}, q_{s}, q_{w})}
\end{align}
\noindent where $\Phi_{\beta}$ - known as the replica symmetric potential - is a concave function depending only on the following scalar parameters:
\begin{align}
q_{s} = \frac{1}{d}||\mat{F}\vec{w}||^2,&& q_{w} = \frac{1}{p}||\vec{w}||^2, && m_{s} = \frac{1}{d}(\mat{F}\vec{w})\cdot \vec{\theta}^{0}\notag
\end{align}
\noindent for $\vec{w}\sim \pi_{\beta}$. In the limit of $p \rightarrow \infty$, this integral concentrates around the extremum of the potential $\Phi_{\beta}^{(0)}$ for any $\beta$. Since the optimisation problem in eq.~\eqref{def-f-fit} is convex, by construction as $\beta\to\infty$ the overlap parameters $(q_{s}^{\star}, q_{w}^{\star}, m_{s}^{\star})$ satisfying this optimisation problem are precisely the ones of eq.~\eqref{eq:optimal_overlaps} corresponding to the solution $\hat{\vec{w}}\in\mathbb{R}^{p}$ of eq.~\eqref{def-f-fit}. 

In summary, the replica method allows to circumvent the hard-to-solve high-dimensional optimisation problem eq.~\eqref{def-f-fit} by directly computing the generalisation error in eq.~\eqref{def-gen} in terms of a simpler scalar optimisation. Doing gradient descent in $\Phi_{\beta}^{(0)}$ and taking $\beta\to\infty$ lead to the saddle-point eqs.~\eqref{eq:saddle_point}.

\subsection{Replicated Gaussian Equivalence}
The backbone of the replica derivation sketched above and detailed in
Sec.~\ref{sec:app:replicas} of the supplementary material is a
central limit theorem type result coined as the \emph{Gaussian
  equivalence theorem} (GET) from \cite{goldt2019modelling}
used in the context of the ``replicated'' Gibbs measure obtained by
taking $r$ copies of (\ref{eq:main:Gibbs_definition}). In this
approach, we need to assume that the ``balance condition''
(\ref{eq:balance}) applies with probability one when the weights $w$
are sampled from the replicated measure. We shall use this assumption
in the following, checking its self-consistency via agreement with
simulations.

 It is interesting to observe that, when applying the GET in the context of our replica calculation, the resulting asymptotic generalisation error stated in Sec.~\ref{sec:main:replica_result} is equivalent to the asymptotic generalisation error of the following linear model: 
\begin{align}
\vec{x}^{\mu} = \kappa_{0} {\bold{1}}+  \kappa_{1} \frac{1}{\sqrt{d}}  \mat{F}^{\top} \vec{c}^{\mu} + \kappa_{\star}  ~\vec{z}^{\mu}\, ,
  \label{model-def-x-GET}
\end{align}
with $\kappa_{0} = \mathbb{E}\left[\sigma(z)\right]$,
$\kappa_{1} \equiv \mathbb{E}\left[z\sigma(z)\right]$,
$\kappa^2_{\star} \equiv
\mathbb{E}\left[\sigma(z)^2\right]-\kappa_{0}^2-\kappa_{1}^2$, and
$\vec{z}^{\mu}\sim {\cal N}(\vec{0},\mat{I}_{p})$. We have for
instance,
$(\kappa_{0}, \kappa_{1},\kappa_{\star}) \approx
\left(0,\frac{2}{\sqrt{3\pi}}, 0.2003 \right)$ for $\sigma = \erf$
and $(\kappa_{0}, \kappa_{1},\kappa_{\star}) = \left(0,
  \sqrt{\frac{2}{\pi}}, \sqrt{1-\frac{2}{\pi}}\right)$ for
$\sigma = \sign$, two cases explored in the next section. This equivalence constitutes a result with an interest in its own, with applicability beyond the scope of the generalised linear task eq.~\eqref{eq-loss} studied here. 

Equation (\ref{model-def-x-GET}) is precisely the mapping obtained by
\cite{mei2019generalisation}, who proved its validity rigorously in
the particular case of the square loss and Gaussian random matrix
$\mat{F}$ using random matrix theory. The same equivalence arises in
the analysis of kernel random matrices \cite{cheng2013,NIPS2017_6857}
and in the study of online learning \cite{goldt2019modelling}. The
replica method thus suggests that the equivalence actually holds in a
much larger class of learning problem, as conjectured as well in
\cite{montanari2019generalisation}, and numerically confirmed in all
our numerical tests. It also potentially allows generalisation of the
analysis in this paper for data coming from a learned generative
adversarial network, along the lines of
\cite{seddik2019kernel,el2020random}.

Fig.~\ref{fig_comparison} illustrates the remarkable agreement between the result
of the generalisation formula, eq.~(\ref{eq:eg_main}) and simulations
both on the data eq.~(\ref{model-def-x}) with $\sigma(x) = {\rm sign}(x)$
non-linearity, and on the Gaussian equivalent data
eq.~(\ref{model-def-x-GET}) where the non-linearity is replaced by
rescaling by a constant plus noise. The agreement is flawless as
implied by the theory in the high-dimensional limit, testifying that
the used system size $d=200$ is sufficiently large for
the asymptotic theory to be relevant. We observed similar good
agreement between the theory and simulation in all the cases we
tested, in particular in all those presented in the following. 
 
\section{Applications of the generalisation formula}

\subsection{Double descent for classification with logistic loss}
\begin{figure}[t]
\begin{center}
\centerline{\includegraphics[width=0.5\columnwidth]{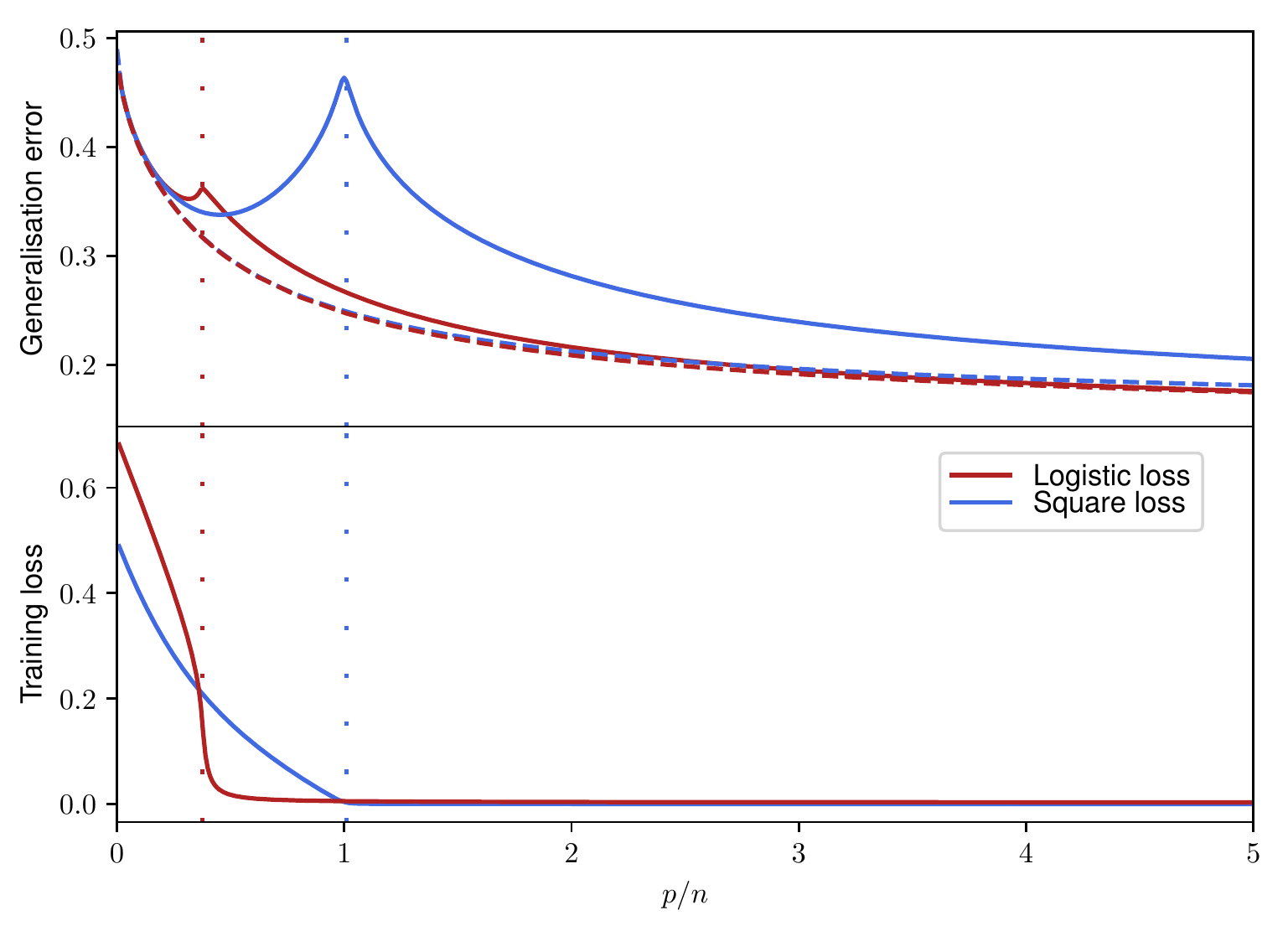}}
\caption{Upper panel: Generalisation error evaluated from eq.~(\ref{eq:eg_main}) plotted against the number of random Gaussian
  features per sample $p/n= 1/\alpha$ and fixed ratio between the
  number of samples and dimension $n/d=\alpha/\gamma=3$ for logistic
  loss (red), square loss (blue). Labels are generated as 
  $y^{\mu}=\sign\left(\vec{c}^{\mu}\cdot\vec{\theta}^{0}\right)$,
  data as $\vec{x}^{\mu}=\sign\left(\mat{F}^{\top}\vec{c}^{\mu}\right)$ and $\hat{f}=\sign$
  for two different values of regularisation  $\lambda$, a small penalty
  $\lambda=10^{-4}$ (full line) and a value of lambda optimised for
  every $p/n$ (dashed line). Lower panel: The training loss corresponding to
  $\lambda=10^{-4}$ is depicted. }
\label{fig:generalisation:logistic:lambdas}
\end{center}
\vskip -0.2in
\end{figure}

\begin{figure}[t]
  \begin{center}
      \vspace{-0.3cm}
  \centerline{\includegraphics[width=0.5\columnwidth]{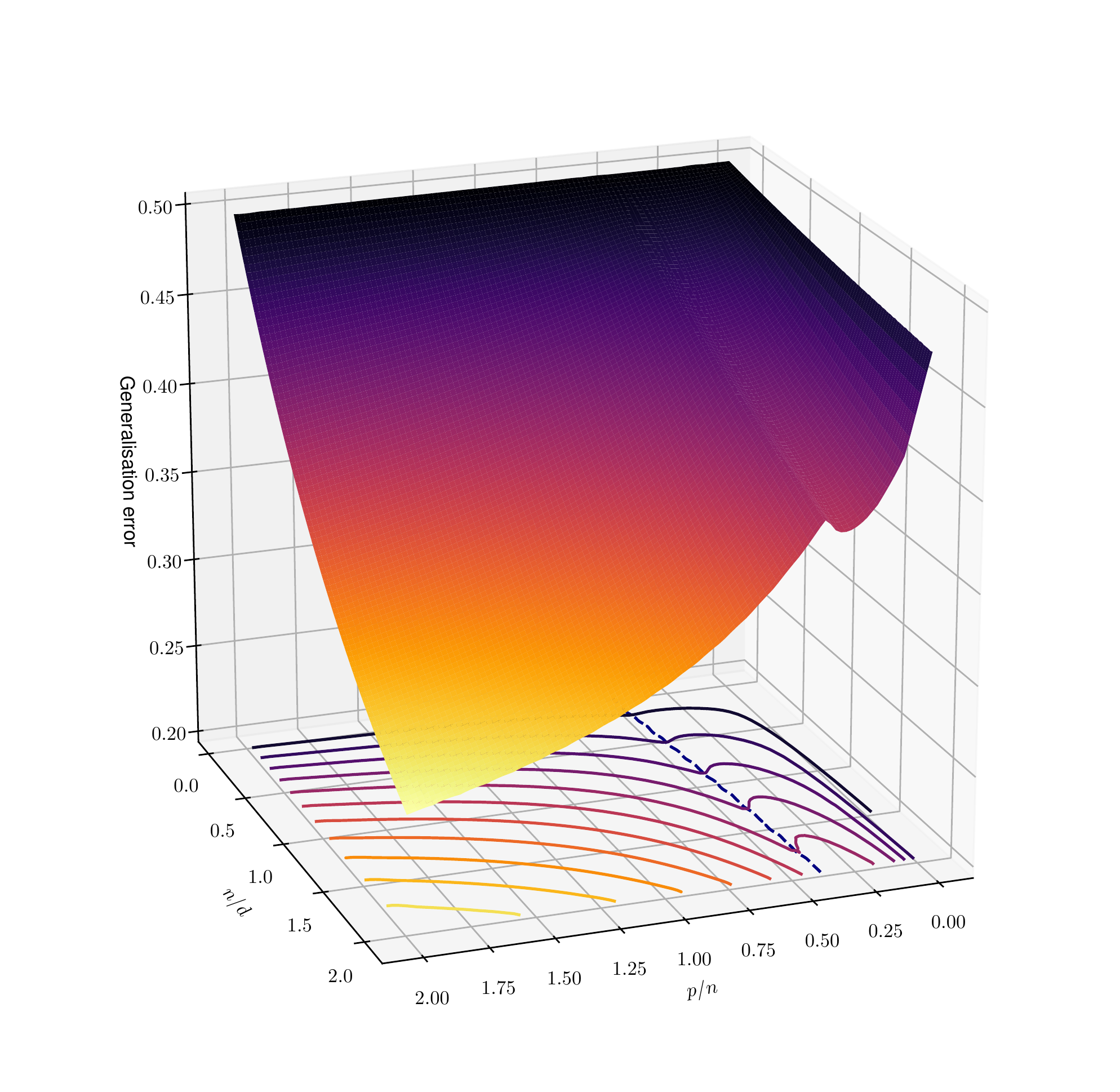}}
  \vspace{-0.8cm}
\caption{Generalisation error of the logistic loss at
  fixed very small regularisation $\lambda = 10^{-4}$, as a function of
  $n/d = \alpha/\gamma$ and $p/n= 1/\alpha$, for random Gaussian
  features. Labels are generated with
  $y^{\mu}=\sign\left(\vec{c}^{\mu}\cdot\vec{\theta}^{0}\right)$,
  the data
  $\vec{x}^{\mu}=\sign\left(\mat{F}^{\top}\vec{c}^{\mu}\right)$ and $\hat{f} = \sign$. The
  interpolation peak happening where data become linearly separable is
  clearly visible here.}
  \label{fig:heatmap:logistic:sign}
\end{center}
\vskip -0.2in
\end{figure}

\begin{figure}[ht]
\begin{center}
\centerline{\includegraphics[width=0.5\columnwidth]{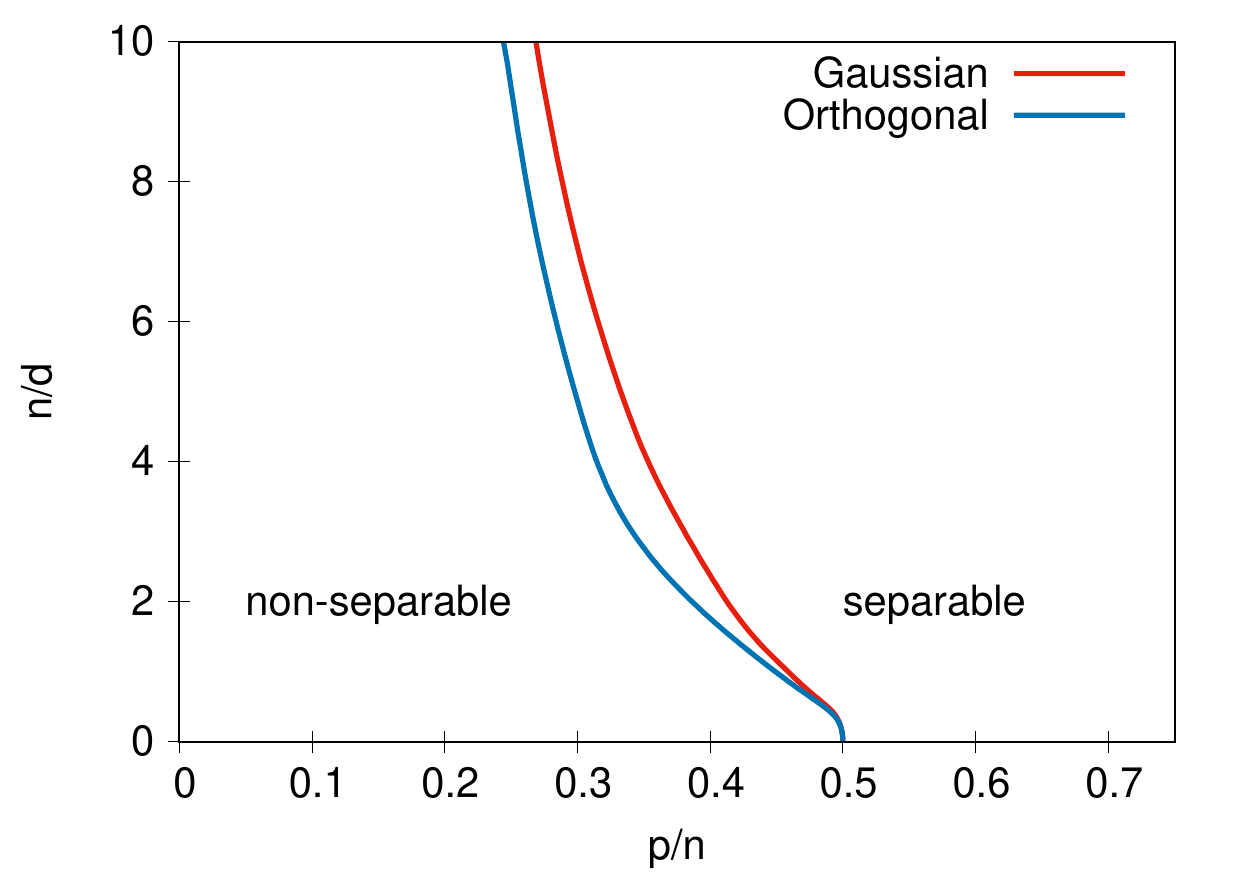}}
\caption{The position of the interpolation peak in
  logistic regression with $\lambda=10^{-4}$, where
  data become linearly separable, as a function of the ratio between
  the number of samples $n$ and the dimension $d$. Labels are generated with
  $y^{\mu}=\sign\left(\vec{c}^{\mu}\cdot\vec{\theta}^{0}\right)$,
  the data
  $\vec{x}^{\mu}=\sign\left(\mat{F}^{\top}\vec{c}^{\mu}\right)$ and $\hat{f} = \sign$. The red line is
  with Gaussian random features, the blue line with orthogonal
  features. We see that for linear separability we need smaller number
  of projections $p$ with orthogonal random features than with
  Gaussian. 
}
\label{fig:generalisation:logistic:phase}
\end{center}
\vskip -0.2in
\end{figure}

Among the surprising observations in modern machine learning is
the fact that one can use learning methods that achieve zero training
error, yet their generalisation error does not deteriorate as more and
more parameters are added into the neural network. The study of such
``interpolators''  have attracted a growing attention over the last few years
\cite{advani2017high, spigler2018jamming,belkin2018reconciling,neal2018modern,hastie2019surprises,mei2019generalisation,geiger2019scaling,nakkiran2019deep},
as it violates basic intuition on the bias-variance trade-off \cite{geman1992neural}. Indeed classical learning theory suggests that generalisation
should first improve then worsen when increasing model complexity,
following a U-shape curve.  Many methods, including neural networks, instead follow a so-called  "double
descent curve"~\cite{belkin2018reconciling} that displays two regimes:
the "classical" U-curve found at low number of parameters is followed at
high number of parameters by an interpolation regime where the
generalisation error decreases monotonically. Consequently 
neural networks do not drastically overfit even when using much more
parameters than data samples \cite{breiman1995reflections}, as
actually observed already in the classical work \cite{geman1992neural}.  Between the two
regimes, a "peak" occurs at the interpolation
threshold~\cite{opper1996statistical,engel2001statistical,
  advani2017high,spigler2018jamming}. It should, however, be noted that existence of this "interpolation" peak is
an independent phenomenon from the lack of overfitting in highly
over-parametrized networks, and indeed in a number of the related works
these two phenomena were observed separately
\cite{opper1996statistical,engel2001statistical,advani2017high,geman1992neural}. Scaling
properties
of the peak and its relation to the jamming phenomena in physics are in
particular studied in \cite{geiger2019scaling}.

Among the simple models that allow to observe this behaviour, random
projections ---that are related to lazy training and kernel
methods--- are arguably the most natural one. The double descent has
been analysed in detail in the present model in the specific case of a square loss
on a regression task with random Gaussian features
\cite{mei2019generalisation}. Our analysis allows to show the
generality and the robustness of the phenomenon to other tasks,
matrices and losses. In Fig.~\ref{fig:generalisation:logistic:lambdas} we compare the
double descent as present in the square loss (blue line) with the one of logistic
loss (red line) for random Gaussian
  features. We plot the value of the generalisation error at
small values of the regularisation $\lambda$ (full line), and for
optimal value of $\lambda$ (dashed line) for a
fixed ratio between the number of samples and the dimension $n/d$ as
a function of the number of random features per sample $p/n$.
We also plot the value of the training error (lower panel) for a
small regularisation value, showing
that the peaks indeed occur when the training loss goes to zero. For the square loss
the peak appears at $1/\alpha\!=\!p/n\!=\!1$ when the system
of $n$ linear equations with $p$ parameters becomes solvable. For the
logistic loss the peak instead appears at a value $1/\alpha^*$ where
the data ${\cal D}$ become linearly separable and hence the logistic
loss can be optimised down to zero. These values $1/\alpha^*$ depends
on the value $n/d$, and this
dependence is plotted in
Fig.~\ref{fig:generalisation:logistic:phase}. For very large dimension
$d$, i.e. $n/d \to 0$ the data matrix $X$ is close to iid random
matrix and hence the $\alpha^*(n/d = 0) = 2$ as famously derived in
classical work by Cover \cite{cover1965geometrical}. For $n/d >0$ the
$\alpha^*$ is growing ($1/\alpha^*$ decreasing) as  correlations make
data easier to linearly separate, similarly as in \cite{candes2020}.

Fig.~\ref{fig:generalisation:logistic:lambdas} also shows that better
error can be achieved with the logistic loss compared to the square
loss, both for small and optimal regularisations, except in a small
region around the logistic interpolation
peak. In the Kernel limit, i.e. $p/n \to \infty$,  the generalization error at optimal regularisation saturates at $\epsilon_g
(p/n \to \infty) \simeq 0.17$ for square loss and at $\epsilon_g
(p/n \to \infty) \simeq 0.16$ for logistic loss. Fig.~\ref{fig:heatmap:logistic:sign} then depicts a 3D plot of
the generalisation error also illustrating the position of the interpolation peak.

\subsection{Random features: Gaussian versus orthogonal}
\begin{figure*}[t!]
\includegraphics[width=0.5\columnwidth]{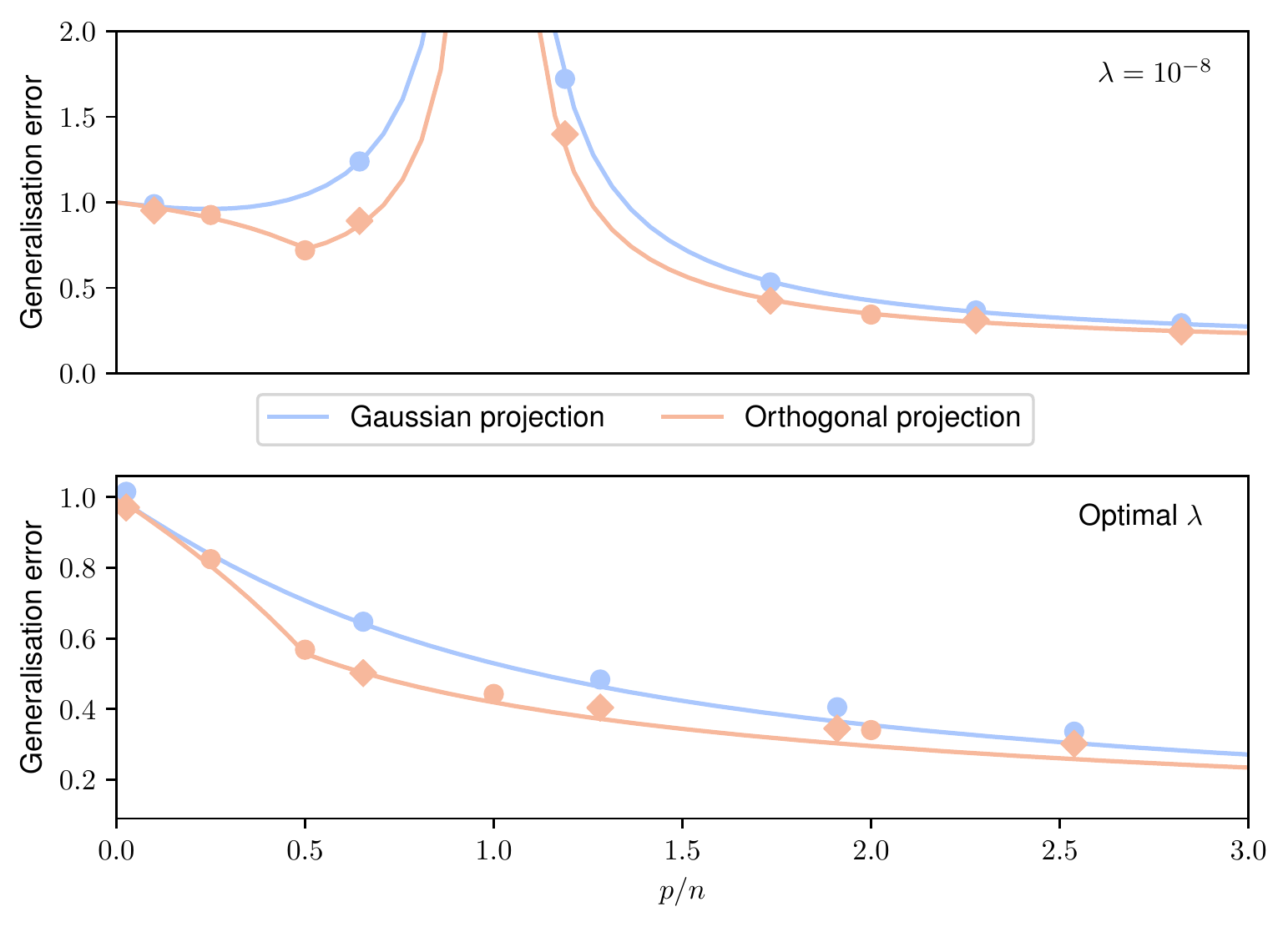}
  \includegraphics[width=0.5\columnwidth]{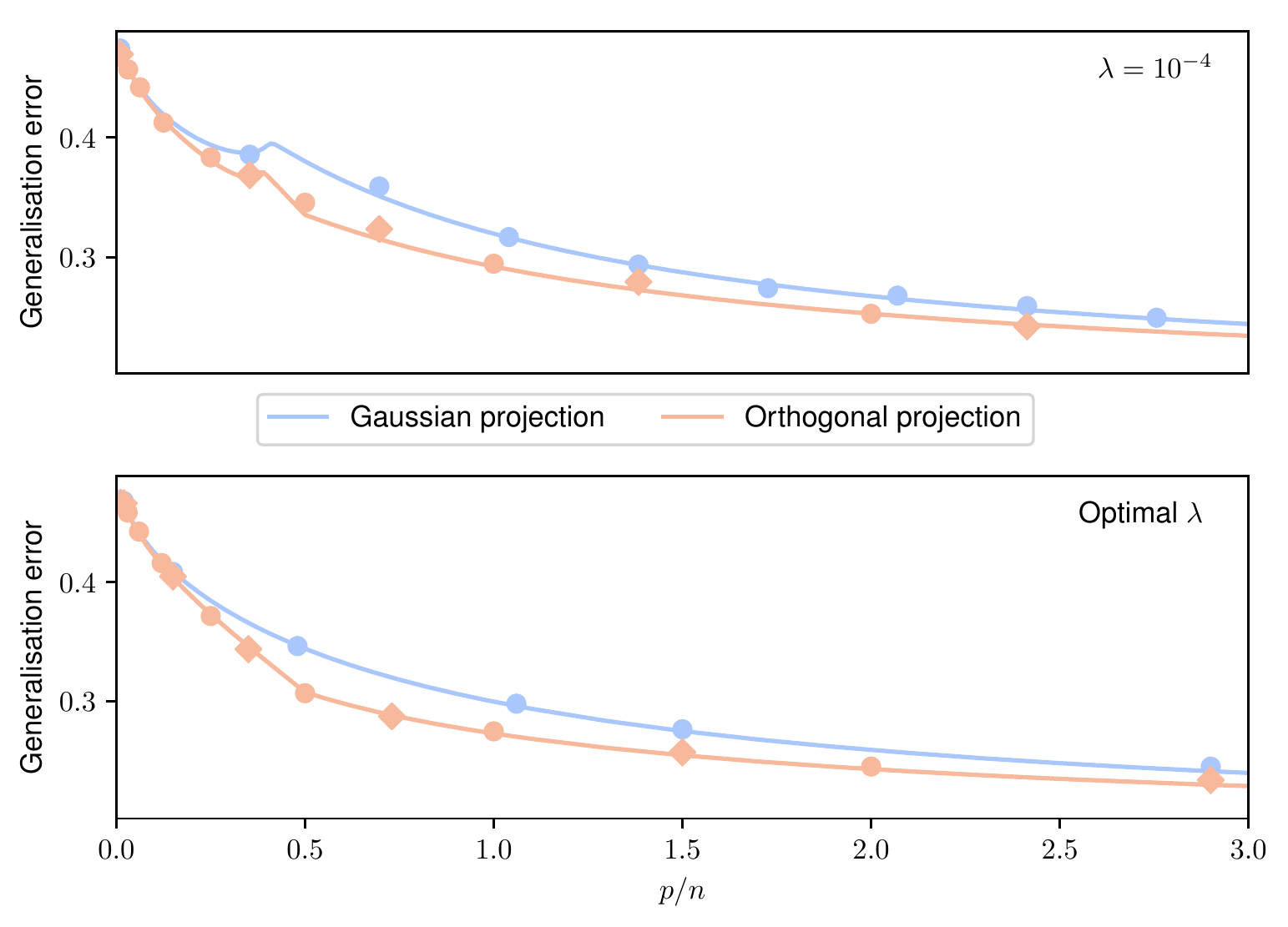}
  \caption{Generalisation error against the number of features per
    sample  $p/n$, for a regression problem (left) and a classification one
    (right). {\bf Left (ridge regression):} We used $n/d=2$ and generated labels as
    $y^{\mu}=\vec{c}^{\mu}\cdot\vec{\theta}^{0}$, data
as  $\vec{x}^{\mu}=\sign\left(\mat{F}^{\top}\vec{c}^{\mu}\right)$ and $\hat{f}(x) = x$. The
  two curves correspond to ridge regression with Gaussian (blue) versus orthogonal (red)
  projection matrix $\mat{F}$ for both $\lambda=10^{-8}$ (top) and
  optimal regularisation $\lambda$ (bottom). {\bf Right (logistic classification)}: We used $n/d=2$
and generated labels as  $y^{\mu}=\sign\left(\vec{c}^{\mu}\cdot\vec{\theta}^{0}\right)$,
  data as
  $\vec{x}^{\mu}=\sign\left(\mat{F}^{\top}\vec{c}^{\mu}\right)$ and $\hat{f}=\sign$. The
  two curves correspond to a logistic classification with again Gaussian (blue) versus orthogonal (red)
  projection matrix $\mat{F}$ for both $\lambda=10^{-4}$ and optimal regularisation
  $\lambda$. In all cases, full lines is the
  theoretical prediction, and points correspond to gradient-descent
  simulations with $d=256$. For the simulations of orthogonally invariant matrices, we results for Hadamard matrices (dots) and DCT Fourier matrices (diamonds).}
\label{fig:generalisation:erf:logistic:optimal}
  \label{fig:generalisation:erf:ridge:optimal}
\end{figure*}

Kernel methods are a very popular class of machine learning techniques,
achieving state-of-the-art performance on a variety of tasks with
theoretical guarantees
\cite{scholkopf2002learning,rudi2017falkon,caponnetto2007optimal}. In
the context of neural network, they are the subject of a renewal of
interest in the context of the Neural Tangent Kernel
\cite{jacot2018neural}. Applying kernel methods to large-scale ``big
data'' problems, however, poses many computational challenges, and
this has motivated a variety of contributions to develop them at
scale, see, e.g.,
\cite{rudi2017falkon,zhang2015divide,saade2016random,ohana2019kernel}. Random
features \cite{NIPS2007_3182} are among the most popular techniques to
do so.

Here, we want to compare
the performance of random projection with respect to
structured ones, and in particular orthogonal random
projections \cite{choromanski2017unreasonable} or deterministic
matrices such as real Fourier (DCT) and Hadamard matrices used in fast
projection methods
\cite{le2013fastfood,andoni2015practical,bojarski2016structured}. Up to normalisation, these matrices have the same spectral density. Since the asymptotic generalisation error only depends on the spectrum of $\mat{F}\mat{F}^{\top}$, all these matrices share the same theoretical prediction when properly normalised, see Fig.~\ref{fig:generalisation:erf:ridge:optimal}. In
our computation, left- and right-orthogonal invariance is parametrised by letting $\mat F = \mat{U}^{\top}\mat{D}\mat{V}$ for $\mat{U}\in\mathbb{R}^{d\times d}$, $\mat{V}\in\mathbb{R}^{p\times p}$ two orthogonal matrices drawn from the Haar measure, and $\mat{D}\in\mathbb{R}^{d\times p}$ a diagonal matrix of rank $\min(d,p)$. In order to compare the results with the Gaussian case, we fix the diagonal entries $d_{k}=\max(\sqrt{\gamma}, 1)$ of $\mat{D}$ such that an arbitrary projected vector has the same norm, on average, to the Gaussian case.

Fig.~\ref{fig:generalisation:erf:ridge:optimal} shows that
random orthogonal embeddings always outperform Gaussian random
projections, in line with empirical observations, and that they allow
to reach the kernel limit with fewer number of projections. Their behaviour is,
however, qualitatively similar to the one of random
i.i.d. projections. We also show in
Fig.~\ref{fig:generalisation:logistic:phase} that orthogonal projections allow to
separate the data more easily than the Gaussian ones, as the phase
transition curve delimiting the linear separability of the logistic loss get shifted to the left.
\begin{figure}[h!]
\begin{center}
\centerline{\includegraphics[width=0.5\columnwidth]{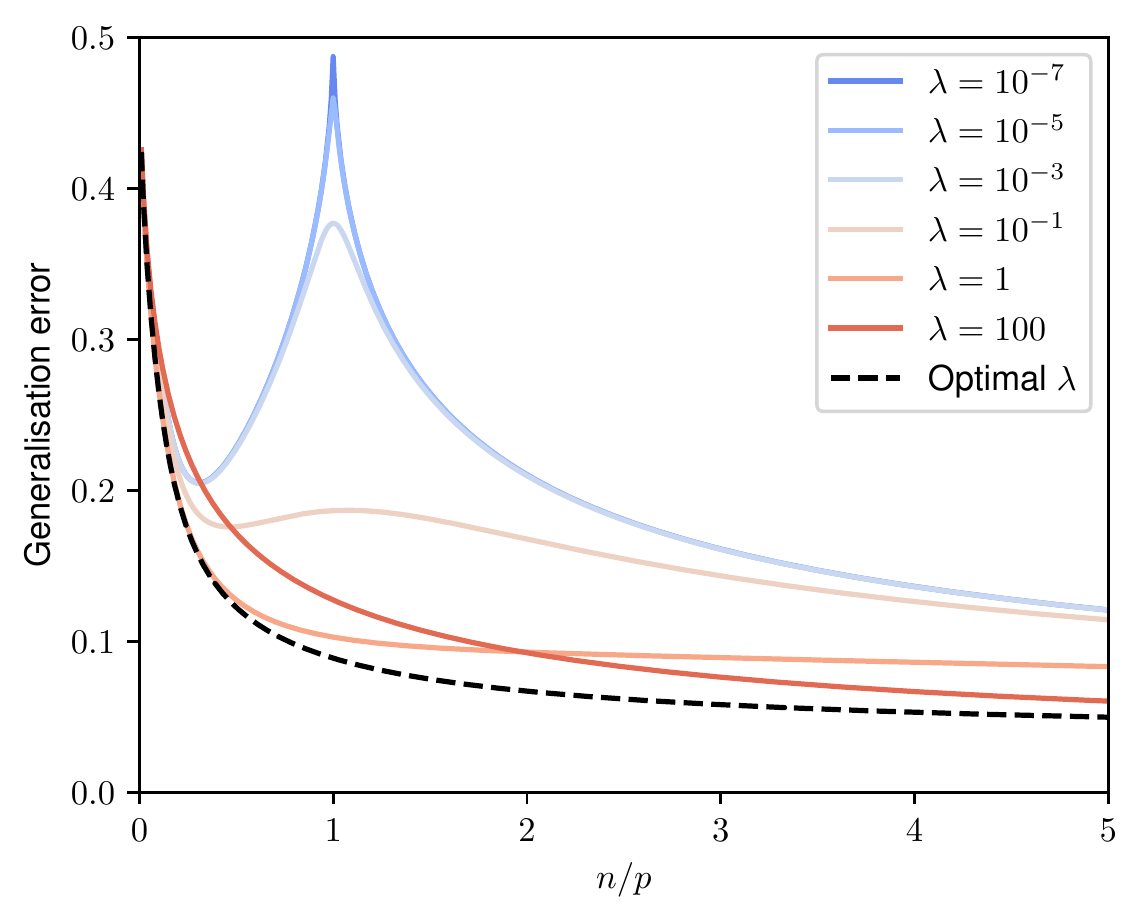}}
\caption{Generalisation error against the number of samples per
  dimension, $\alpha= n/p$, and
  fixed ratio between the latent and data dimension, $d/p=0.1$, for a
  classification task with square loss on labels generated as
  $y^{\mu}=\sign\left(\vec{c}^{\mu}\cdot\vec{\theta}^{0}\right)$,
  data $\vec{x}^{\mu}=\erf\left(\mat{F}^{\top} \vec{c}^{\mu}\right)$ and $\hat{f}=\sign$,
  for different values of the regularisation $\lambda$ (full lines), including the
  optimal regularisation value (dashed).}
\label{fig:generalisation:logistic:HMM}
\end{center}
\vskip -0.2in
\end{figure}

\begin{figure}[h!]
\begin{center}
\centerline{\includegraphics[width=0.5\columnwidth]{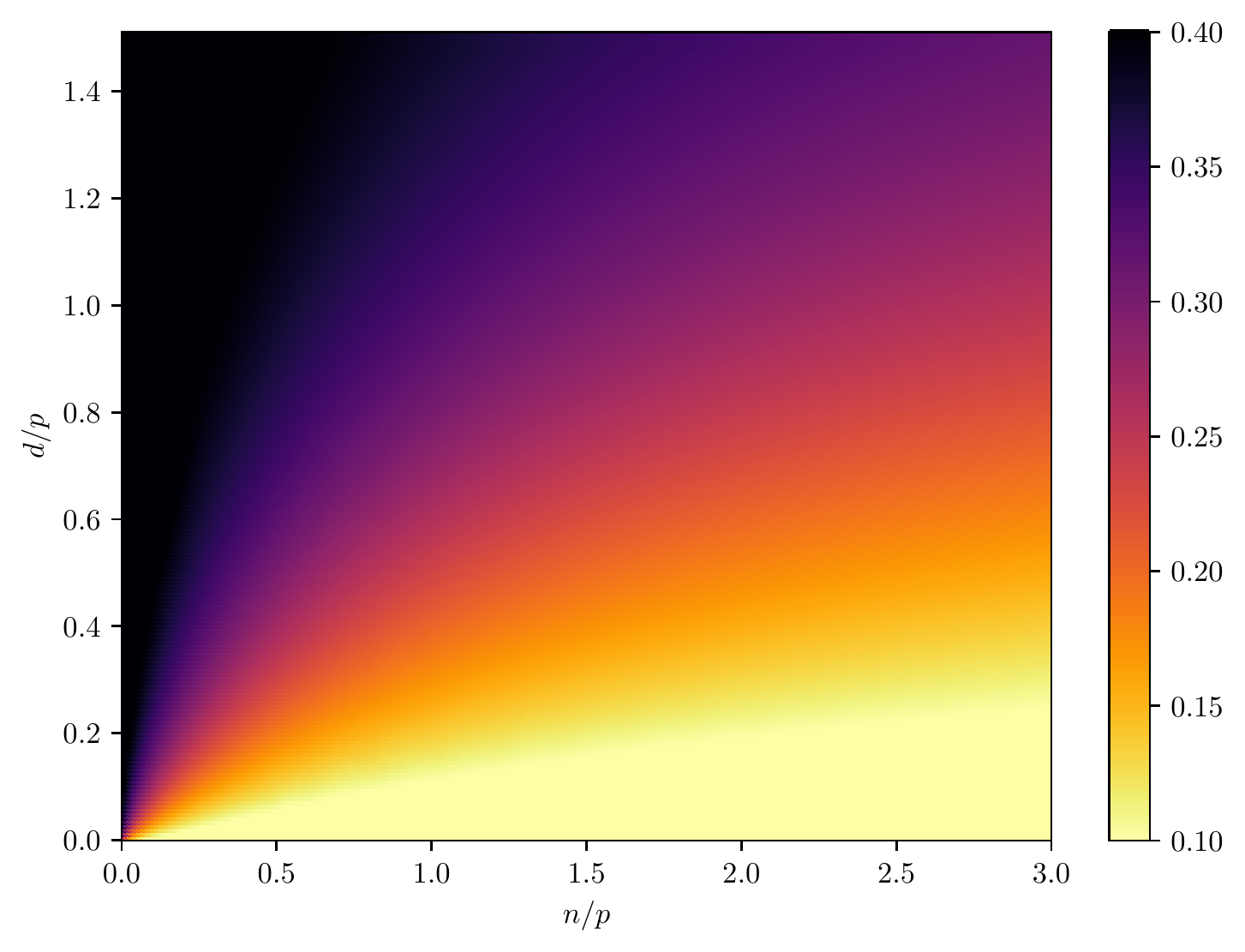}}
\caption{Heat-map of the generalisation errors as a function of the
  number of 
  samples per data dimension $n/p$ against the ratio of the latent and
  data dimension $d/p$, for a classification task with square loss on labels
  $y^{\mu}=\sign\left(\vec{c}^{\mu}\cdot\vec{\theta}^{0}\right)$ and
  data $\vec{x}^{\mu}=\erf\left(\mat{F}^{\top} \vec{c}^{\mu}\right)$
  for the optimal values of the regularisation $\lambda$.}
\label{fig:generalisation:logistic:HMM_heat}
\end{center}
\vskip -0.3in
\end{figure}

\subsection{The hidden manifold model phase diagram}
In this subsection we consider the hidden manifold model where
$p$-dimensional $x$ data lie on a $d$-dimensional manifold, we have mainly
in mind $d < p$. The labels $y$ are generated using the coordinates on
the manifold, eq.~(\ref{model-def-y}). 

In Fig.~\ref{fig:generalisation:logistic:HMM} we plot the
generalisation error of classification with the square loss
for
various values of the
regularisation $\lambda$. We fix the ratio between the dimension of
the sub-manifold and the dimensionality of the input data to $d/p=0.1$
and plot the learning curve, i.e. the error as a function of the
number of samples per dimension. Depending on the value of the
regularisation, we observe that the interpolation peak, which is at $\alpha =1$ at
very small regularisation (here
the over-parametrised regime is on the left hand side), 
decreases for larger regularisation~$\lambda$. A similar behaviour has
been observed for other models in the past, see e.g.~\cite{opper1996statistical}.
Finally Fig.~\ref{fig:generalisation:logistic:HMM} depicts the error for
optimised regularisation parameter in the black dashed
line. For large number of samples we observe the generalisation error
at optimal regularisation to saturate in this case at $\epsilon_g
(\alpha \to \infty) \to 0.0325$.
A challenge for future work is to
see whether better performance can be achieved on this model by including hidden variables into the neural network.

Fig.~\ref{fig:generalisation:logistic:HMM_heat} then shows the
generalisation error for the optimised regularisation $\lambda$ with
square loss as a function of the ratio between the latent and the data
dimensions $d/p$. In the limit $d/p \gg 1$ the data matrix becomes close to
a random iid matrix and the labels are effectively random, thus only bad
generalisation can be reached. Interestingly, as $d/p$ decreases to
small values even the simple classification with regularised square
loss is able to ``disentangle'' the hidden manifold structure in the
data and to reach a rather low generalisation error. The figure
quantifies how the error deteriorates
when the ratio between the
two dimensions $d/p$ increases. Rather remarkably, for a low $d/p$ a good
generalisation error is achieved even
in the over-parametrised regime, where the dimension is
larger than the number of samples, $p > n$. In a sense, the square loss
linear classification is able to locate the low-dimensional subspace and find
good generalisation even in the over-parametrised regime as long as
roughly $d \lesssim n$. The observed results are in qualitative agreement with
the results of learning with stochastic gradient descent in
\cite{goldt2019modelling} where for very low $d/p$ good
generalisation error was observed in the hidden manifold model, but a
rather bad one for $d/p=0.5$.

\section*{Acknowledgements}
This work is supported by the ERC under the European Union's Horizon
2020 Research and Innovation Program 714608-SMiLe, as well as by the
French Agence Nationale de la Recherche under grant
ANR-17-CE23-0023-01 PAIL and ANR-19-P3IA-0001 PRAIRIE. We also
acknowledge support from the chaire CFM-ENS  "Science des donn\'ees''.
We thank Google Cloud for providing us access to their platform through the Research Credits Application program. BL was partially financed by the Coordena\c{c}\~ao de Aperfei\c{c}oamento de Pessoal de N\'ivel Superior - Brasil (CAPES) - Finance Code 001.

\newpage
\appendix
\section*{Appendix}
\section{Definitions and notations}
\label{sec:app:definitions_notations}
In this section we recall the models introduced in the main body of the article, and introduce the notations used throughout the Supplementary Material.

\subsection{The dataset}
\label{sec:app:data}
In this work we study a series of regression and classification tasks for a dataset $\{\vec{x}^{\mu}, y^{\mu}\}_{\mu=1}^{n}$ with labels $y^{\mu}\in\mathbb{R}$ sampled identically from a \emph{generalised linear model}:
\begin{align}
	y^{\mu} \sim P^{0}_{y}\left(y^{\mu}\Big|\frac{\vec{c}^{\mu}\cdot \vec{\theta}^{0}}{\sqrt{d}}\right),
	\label{eq:app:rule}
\end{align}
\noindent  where the output-channel $P^{0}_{y}\left( \cdot \right) $ is defined as:
\begin{equation}
P^{0}_{y}\left(y^{\mu}\Big|\frac{\vec{c}^{\mu}\cdot \vec{\theta}^{0}}{\sqrt{d}}\right) = \int d \xi^{\mu} P\left( \xi^{\mu}\right) \delta \left( y^{\mu} - f^0\left( \frac{\vec{c}^{\mu}\cdot \vec{\theta}^{0}}{\sqrt{d}}; \xi^{\mu}\right)     \right)   
\end{equation}
\noindent for some noise $\xi^{\mu}$ and for data points $\vec{x}^{\mu}\in\mathbb{R}^{p}$ given by:
\begin{align}
\vec{x}^{\mu} = \sigma\left(\frac{1}{\sqrt{d}}\sum\limits_{\rho=1}^{d}c^{\mu}_{\rho}\vec{f}_{\rho}\right).	
\end{align}
The vectors $\vec{c}^{\mu}\in\mathbb{R}^{d}$ is assumed to be identically drawn from $\mathcal{N}(0,\mat{I}_{d})$, and $\vec{\theta}^{0}\in\mathbb{R}^{d}$ from a separable distribution $P_{\theta}$. The family of vectors $\vec{f}_{\rho}\in\mathbb{R}^{p}$ and the scalar function $\sigma:\mathbb{R}\to\mathbb{R}$ can be arbitrary. 

Although our results are valid for the general model introduced above, the two examples we will be exploring in this work are the noisy linear channel (for regression tasks) and the deterministic sign channel (for classification tasks):
\begin{align}
y^{\mu} = \frac{\vec{c}^{\mu}\cdot \vec{\theta}^{0}}{\sqrt{d}}+\sqrt{\Delta}~\xi^{\mu}\qquad &\Leftrightarrow\qquad P^{0}_{y}\left(\vec{y}\Big|\frac{\vec{c}^{\mu}\cdot \vec{\theta}^{0}}{\sqrt{d}}\right)= \prod\limits_{\mu=1}^{n}\mathcal{N}\left(y^{\mu}; \frac{\vec{c}^{\mu}\cdot \vec{\theta}^{0}}{\sqrt{d}}, \Delta\right)\\
y^{\mu} = \sign\left(\frac{\vec{c}^{\mu}\cdot \vec{\theta}^{0}}{\sqrt{d}}\right) \qquad &\Leftrightarrow\qquad P^{0}_{y}\left(\vec{y}\Big|\frac{\vec{c}^{\mu}\cdot \vec{\theta}^{0}}{\sqrt{d}}\right) = \prod\limits_{\mu=1}^{n}\delta\left(y^{\mu}-\sign\left(\frac{\vec{c}^{\mu}\cdot \vec{\theta}^{0}}{\sqrt{d}}\right)\right)
\end{align}
\noindent where $\xi^{\mu}\sim\mathcal{N}(0,1)$ and $\Delta > 0$. 

This dataset can be regarded from two different perspectives.
\paragraph{Hidden manifold model:} The dataset $\{\vec{x}^{\mu}, y^{\mu}\}_{\mu=1,\cdots,n}$ is precisely the \emph{hidden manifold model} introduced in \cite{goldt2019modelling} to study the dynamics of online learning in a synthetic but structured dataset. From this perspective, although $\vec{x}^{\mu}$ lives in a $p$ dimensional space, it is parametrised by a latent $d<p$-dimensional subspace spanned by the basis $\{\vec{f}_{\rho}\}_{\rho=1,\cdots, d}$ which is "hidden" by the application of a scalar nonlinear function $\sigma$ acting component-wise. The labels $\vec{y}^{\mu}$ are then drawn from a generalised linear rule defined on the latent $d$-dimensional space. 

\paragraph{Random features model:} The dataset $\{\vec{x}^{\mu}, y^{\mu}\}_{\mu=1,\cdots,n}$ is tightly related to the Random Features model studied in \cite{NIPS2007_3182} as a random approximation for kernel ridge regression. In this perspective, $\vec{c}^{\mu}\in\mathbb{R}^{d}$ is regarded as a collection of $d$-dimensional data points which are projected by a random feature matrix $\mat{F} = (\vec{f}_{\rho})_{\rho=1}^{p}\in\mathbb{R}^{d\times p}$ into a higher dimensional space, followed by a non-linearity $\sigma$. In the limit of infinite number of features $d,p\to\infty$ with fixed ratio $d/p$, performing ridge regression of $\vec{x}^{\mu}$ is equivalent to kernel ridge regression with a limiting kernel depending on the distribution of the feature matrix $\mat{F}$ and on the non-linearity $\sigma$.

\subsection{The task}
In this work, we study the problem of learning the rule from eq.~\eqref{eq:app:rule} from the dataset $\{(\vec{x}^{\mu}, y^{\mu})\}_{\mu=1,\cdots,n}$ introduced above with a \emph{generalised linear model}:
\begin{align}
\hat{y}^{\mu} = \hat{f}\left(\vec{x}^{\mu}\cdot \hat{\vec{w}}\right)
\end{align}
\noindent where the weights $\vec{w}\in\mathbb{R}^{p}$ are learned by minimising a loss function with a ridge regularisation term:
\begin{align}
	\hat{\vec{w}} = \underset{\vec{w}}{\min}~ \left[\sum\limits_{\mu=1}^{n}\ell(y^{\mu}, \vec{x}^{\mu} \cdot \vec{w} )+\frac{\lambda}{2}||\vec{w}||_{2}^2\right] \ .
	\label{eq:app:minimisation}
\end{align}
\noindent for $\lambda > 0$.

It is worth stressing that our results hold for general $\ell$, $\hat{f}$ and $f^0$ - including non-convex loss functions. However, for the purpose of the applications explored in this manuscript, we will be mostly interested in the cases $\hat{f}(x) = f^0(x) = x$ for regression and $\hat{f}(x) = f^0(x)=\sign(x)$ for classification, and we will focus on the following two loss functions:
\begin{align}
\ell(y^{\mu},\vec{x}^{\mu} \cdot \vec{w}) = \begin{cases}
 \frac{1}{2}(y^{\mu} - \vec{x}^{\mu} \cdot \vec{w})^2, & \text{ square loss}\\
 \log\left(1+e^{-y^{\mu} \left(\vec{x}^{\mu} \cdot \vec{w} \right) }\right), & \text{ logistic loss }
 \end{cases}
\end{align}
Note that these loss functions are strictly convex. Therefore, for these losses, the regularised optimisation problem in (\ref{eq:app:minimisation}) has a unique solution.

Given a new pair $(\vec{x}^{\new}, y^{\new})$ drawn independently from the same distribution as $\{(\vec{x}^{\mu}, y^{\mu})\}_{\mu=1}^{n}$, we define the success of our fit through the generalisation error, defined as:
\begin{align}
	\epsilon_{g} = \frac{1}{4^{k}}\mathbb{E}_{\vec{x}^{\new}, y^{\new}}\left(y^{\new}-\hat{y}^{\new}\right)^2
	\label{eq:app:generalisation_error}
\end{align}
\noindent where $\hat{y}^{\new} = \hat{f}(\vec{x}^{\new}\cdot \hat{\vec{w}})$, and for convenience we choose $k=0$ for the regression tasks and $k=1$ for the classification task, such that the generalisation error in this case counts misclassification. Note that for a classification problem, the generalisation error is just one minus the classification error.

Similarly, we define the \emph{training loss} on the dataset $\{\vec{x}^{\mu}, y^{\mu}\}_{\mu=1}^{n}$ as:
\begin{align}
\epsilon_{t} = \frac{1}{n} \mathbb{E}_{\lbrace\vec{x}^{\mu},y^{\mu}\rbrace}\left[ \sum_{\mu = 1}^{n} \ell\left( y^{\mu}, \vec{x}^{\mu} \cdot \vec{\hat{w}}\right) + \frac{\lambda}{2} \Vert \vec{\hat{w}} \Vert_2^2 \right].
\end{align}
Finally, all the results of this manuscript are derived in the \emph{high-dimensional limit}, also known as \emph{thermodynamic limit} in the physics literature, in which we take $p,d,n\to\infty$ while keeping the ratios $\alpha = n/p$, $\gamma = d/p$ fixed.

\section{Gaussian equivalence theorem}
\label{sec:app:get}
In this section we introduce the \emph{replicated Gaussian equivalence} (rGE), a central result we will need for our replica calculation of the generalisation error in Sec.~\ref{sec:main:replica_result} of the main body. The rGET is a stronger version of the Gaussian equivalence theorem (GET) that was introduced and proved in \cite{goldt2019modelling}. Previously, particular cases of the GET were derived in the context of random matrix theory \cite{hachem2007, cheng2013, Fan2019, NIPS2017_6857}. The gaussian equivalence has also been stated and used in \cite{mei2019generalisation, montanari2019generalisation}. 

\subsection{Gaussian equivalence theorem}
Let $\mat{F}\in\mathbb{R}^{d\times p}$ be a fixed matrix, $\vec{w}^{a}\in\mathbb{R}^{p}$, $1\leq a\leq r$ be a family of vectors, $\vec{\theta}^0\in\mathbb{R}^{d}$ be a fixed vector and $\sigma:\mathbb{R}\to\mathbb{R}$ be a scalar function acting component-wise on vectors.

Let $\vec{c}\in\mathbb{R}^{d}$ be a Gaussian vector $\mathcal{N}(0,\mat{I}_{d})$. The GET is a statement about the (joint) statistics of the following $r+1$  random variables
\begin{align}
\lambda^{a} = \frac{1}{\sqrt{p}}\vec{w}^{a}\cdot \sigma(\vec{u})\in\mathbb{R}, && \nu = \frac{1}{\sqrt{d}}\vec{c}\cdot \vec{\theta}^0\in\mathbb{R}.
\end{align}
\noindent in the asymptotic limit where $d,p\to\infty$ with fixed $p/d$ and fixed $r$. For simplicity, assume that $\sigma(x)=-\sigma(-x)$ is an odd function. Further, suppose that in the previously introduced limit the following two balance conditions hold:

Condition 1:
\begin{align}
 \frac{1}{\sqrt{d}}\sum_{\rho = 1}^d F_{i\rho}F_{j\rho}=O(1),
\label{eq:app:balance1}
\end{align}
\noindent for any $\rho$.

Condition 2:
\begin{align}
S^{a_1,...,a_k}_{\rho_1,...,\rho_q}=\frac{1}{\sqrt{p}}\sum\limits_{i=1}^{p}w_{i}^{a_{1}}w_{i}^{a_{2}}\cdots w_{i}^{a_{k}} F_{i\rho_{1}}F_{i\rho_{2}}\cdots F_{i\rho_{q}} = O(1), 
\label{eq:app:ballance}
\end{align}
\noindent for any integers  $k \ge 0$, $q >0$, for any choice of indices $\rho_{1}, \rho_{2}, \cdots, \rho_{q}\in\{1,\cdots, d\}$ all distinct from each other, and for any choice of indices $a_{1}, a_{2}, \cdots, a_{k}\in\{1,\cdots, r\}$. Under the aforementioned conditions, the following theorem holds:
\begin{theorem}
\label{thm:GET}
In the limit $d,p\to\infty$ with fixed $p/d$, the random variables $\{\lambda^{a},\nu\}$ are jointly normal, with zero mean and covariances:
\begin{align}
	 \mathbb{E}\left[\lambda^{a}\lambda^{b}\right]= \frac{\kappa_{\star}^2}{p}\vec{w}^{a}\cdot\vec{w}^{b} +\frac{\kappa_{1}^{2}}{d}\vec{s}^{a}\cdot\vec{s}^{b}, &&
\mathbb{E}\left[\nu^2\right] = \frac{1}{d}||\vec{\theta}^0||^2\notag
\end{align}
\begin{align}
\mathbb{E}\left[\lambda^{a} \nu\right] = \frac{\kappa_{1}}{d}\vec{s}^{a}\cdot \vec{\theta}^0
\end{align}
\noindent where:
\begin{align}
\vec{s}^{a} = \frac{1}{\sqrt{p}}\mat{F}\vec{w}^{a}\in\mathbb{R}^{d}, \qquad a=1,\cdots, r
\end{align}
\noindent and
\begin{align}
\kappa_{0} = \mathbb{E}_{z}\left[\sigma(z)\right], &&	 \kappa_{1} = \mathbb{E}_{z}\left[z\sigma(z)\right], && \kappa_{\star}^2 = \mathbb{E}_{z}\left[\sigma(z)^2\right]-\kappa_{0}^2-\kappa_{1}^2
\label{eq:app:kappas}
\end{align}
\noindent where $z\sim\mathcal{N}(0,1)$.
\end{theorem}

\subsection{Replicated Gaussian equivalence}
\label{sec:app:rget}
Note that the GET holds for a fixed family $\{\vec{w}^{a}\}_{a=1}^{r}$ and matrix $\mat{F}\in\mathbb{R}^{d\times p}$ satisfying the balance condition from eq.~\eqref{eq:app:ballance}. In the replica setting, we will need to apply the GET under an average over $r$ samples (refered here as \emph{replicas}) of the Gibbs distribution $\mu_{\beta}$, introduced in eq.~\ref{eq:main:Gibbs_definition} on the main. We therefore shall require the assumption that the balance condition eq.~\eqref{eq:app:ballance} holds for any sample of $\mu_{\beta}$. We refer to this stronger version of the GET as the \emph{replicated Gaussian equivalence} (rGE).
Although proving this result is out of the scope of the present work, we check its self-consistency extensively with numerical simulations.

\section{Replica analysis}
\label{sec:app:replicas}
In this section we give a full derivation of the result in Sec. \ref{main-sec:main:replica_result} in the main manuscript for the generalisation error of the problem defined in Sec. \ref{sec:app:definitions_notations}. Our derivation follows from a Gibbs formulation of the optimisation problem in eq. ~\eqref{eq:app:minimisation} followed by a replica analysis inspired by the toolbox of the statistical physics of disordered systems.

\subsection{Gibbs formulation of problem}
Given the dataset
$\{\vec{x}^{\mu}, y^{\mu}\}_{\mu=1}^{n}$ defined in Section \ref{sec:app:data}, we define the following Gibbs measure over $\mathbb{R}^{p}$:
\begin{align}
  \mu_{\beta}(\vec{w}|\{\vec{x}^{\mu},y^{\mu}\})	 = \frac{1}{\mathcal{Z}_{\beta}}e^{-\beta\left[\sum\limits_{\mu=1}^{n}\ell\left(y^{\mu}, \vec{x}^{\mu}\cdot \vec{w}\right) + \frac{\lambda}{2}||\vec{w}||_{2}^2\right]} = \frac{1}{\mathcal{Z}_{\beta}}\underbrace{\prod\limits_{\mu=1}^{n}e^{-\beta\ell\left(y^{\mu}, \vec{x}^{\mu}\cdot \vec{w}\right)}}_{\equiv P_{y}(\vec{y}|\vec{w}\cdot \vec{x}^{\mu})}\underbrace{\prod\limits_{i=1}^{p}e^{-\frac{\beta\lambda}{2} w_{i}^{2}}}_{\equiv P_{w}(\vec{w})}
  \label{eq:Gibbs_definition}
\end{align}
\noindent for $\beta>0$.
When $\beta\to\infty$, the Gibbs measure peaks at the solution of the optimisation problem in eq.~\eqref{eq:app:minimisation} - which, in the particular case of a strictly convex loss, is unique. Note that in the second equality we defined the factorised distributions $P_{y}$ and $P_{w}$, showing that $\mu_{\beta}$ can be interpreted as a posterior distribution of $\vec{w}$ given the dataset $\{\vec{x}^{\mu},y^{\mu}\}$, with  $P_{y}$ and $P_{w}$ being the likelihood and prior distributions respectively. 

An exact calculation of $\mu_{\beta}$ is intractable for large values of $n, p$ and $d$. However, the interest in $\mu_{\beta}$ is that in the limit $n,p,d\to\infty$ with $d/p$ and $n/p$ fixed, the free energy density associated to the Gibbs measure:
\begin{align}
	f_{\beta} = -\lim\limits_{p\to\infty}\frac{1}{p}\mathbb{E}_{\lbrace\vec{x}^{\mu},y^{\mu}\rbrace}\log\mathcal{Z}_{\beta}
	\label{eq:app:def_freeen}
\end{align}
\noindent can be computed exactly
using the replica method,
and at $\beta\to\infty$ give us the optimal overlaps:
\begin{align}
q_{w} = \frac{1}{p}\mathbb{E}||\vec{\hat{w}}||^2	 && q_{x} = \frac{1}{d}\mathbb{E}||\mat{F}\vec{\hat{w}}||^2	 && m_{x} = \frac{1}{d}\mathbb{E}\left[\vec{\theta}^{0}\cdot \mat{F}\hat{\vec{w}}\right]
\end{align}
\noindent that - as we will see - fully characterise the generalisation error defined in eq.~\eqref{eq:app:generalisation_error}.

\subsection{Replica computation of the free energy density}
The replica calculation of $f_\beta$ is based on a large deviation principle for the free energy density. Let
\begin{align}
	f_{\beta}(\{\vec{x}^{\mu},y^{\mu}\}) = -\frac{1}{p} \log\mathcal{Z}_{\beta}
	\label{eq:app:def_freeen_onesample}
\end{align}
\noindent be the free energy density for one given sample of the problem, i.e. a fixed dataset  $\{\vec{x}^{\mu}, y^{\mu}\}_{\mu=1}^{n}$. We assume that the distribution $P(f)$ of the free energy density, seen as a random variable over different samples of the problem, satisfies a large deviation principle, in the sense that, in the thermodynamic limit:
\begin{align}
  P(f)\simeq e^{p \Phi(f)}\ ,
\label{eq:ldp}
\end{align}
with $\Phi$ a concave function reaching its maximum at the free energy density $f=f_\beta$, with $\Phi(f_\beta)=0$.
This hypothesis includes the notion of  \emph{self-averageness} which states that the free-energy density is the same for almost all samples in the thermodynamic limit.

The value of $f_\beta$ can be computed by computing the \emph{replicated partition function}
\begin{align}
\mathbb{E}_{\lbrace\vec{x}^{\mu},y^{\mu}\rbrace}\mathcal{Z}_{\beta}^{r}= \int df \; e^{p[\Phi(f)-r f]}\ ,
\end{align}
and taking the limit
\begin{align}
f_\beta= \lim\limits_{r\to 0^{+}}\frac{\dd}{\dd r}\lim\limits_{p\to\infty}\left[ -\frac{1}{p} \left(\mathbb{E}_{\lbrace\vec{x}^{\mu},y^{\mu}\rbrace}\mathcal{Z}_{\beta}^{r}\right)\right]
 \end{align} 
Although this procedure is not fully rigorous, experience from the statistical physics of disordered systems shows that it gives exact results, and in fact the resulting expression can be verified to match the numerical simulations. 

Using the replica method we need to evaluate:
\begin{align}
  \mathbb{E}_{\lbrace\vec{x}^{\mu},y^{\mu}\rbrace}\mathcal{Z}_{\beta}^{r} &= \int\dd\vec{\theta}^{0}~P_{\theta}(\vec{\theta}^{0})\int\prod\limits_{a=1}^{r}\dd\vec{w}~P_{w}\left(\vec{w}^{a}\right)\times\notag
  \\
                                          &\qquad\times\prod\limits_{\mu=1}^{n}\int\dd y^{\mu}~\underbrace{\mathbb{E}_{\vec{c}^{\mu}}\left[P_{y}^{0}\left(y^{\mu}\big|\frac{\vec{c}^{\mu}\cdot\vec{\theta}^{0}}{\sqrt{d}}\right)\prod\limits_{a=1}^{r}P_{y}\left(y^{\mu}|\vec{w}^{a}\cdot \sigma\left(\frac{1}{\sqrt{d}}\mat{F}^{\top}\vec{c}^{\mu}\right)  \right)\right]}_{\text{(I)}}
                             \label{eq:Zr-av}
\end{align}
\noindent
where $P_w$ and $P_y$ have been defined in (\ref{eq:Gibbs_definition}).
In order to compute this quantity, we introduce, for each point $\mu$ in the database, the $r+1$ variables
\begin{align}
  \nu_{\mu}&=\frac{1}{\sqrt{d}}\vec{c}^{\mu}\cdot\vec{\theta}^{0}\ ,\\
  \lambda_{\mu}^{a}&=\vec{w}^{a}\cdot \sigma\left(\frac{1}{\sqrt{d}}\mat{F}^{\top}\vec{c}^{\mu}\right)\ .
\end{align}
\noindent
Choosing $\vec{c}^{\mu}$ at random induces a joint distribution $P(\nu_{\mu},\lambda^{a}_{\mu})$. In the thermodynamic limit $p,d\to\infty$ with fixed $p/n$, and for matrices $\mat{F}$ satisfying the balance condition in eq.~\eqref{eq:app:ballance}, the \emph{replicated Gaussian equivalence} introduced in Section~\ref{sec:app:rget} tells us that, for a given $\mu$, the $r+1$ variables $\{\nu_{\mu}, \lambda^{a}_{\mu}\}_{a=1}^{r}$ are Gaussian random values with zero mean and covariance given by:
\begin{align}
	\Sigma^{ab} = 
	\begin{pmatrix}
		\rho & M^{a} \\
		M^{a} & Q^{ab}
	\end{pmatrix}\in\mathbb{R}^{(r+1)\times (r+1)}
	\label{eq:app:correlation}
\end{align}
\noindent The elements of the covariance matrix $M^a$ and $Q^{ab}$ are the rescaled version of the so-called \emph{overlap parameters}:
\begin{align}
\rho~ = \frac{1}{d}||\vec{\theta}^{0}||^2, && m_{s}^{a} = \frac{1}{d}\vec{s}^{a}\cdot \vec{\theta}^{0}, && q_{s}^{ab} = \frac{1}{d}\vec{s}^{a}\cdot \vec{s}^{b}, && q_{w}^{ab} = \frac{1}{p}\vec{w}^{a}\cdot \vec{w}^{b}, 	
\end{align}
\noindent where $\vec{s}^{a} = \frac{1}{\sqrt{p}}\mat{F}\vec{w}^{a}$. They are thus given by:
\begin{align}
M^{a} = \kappa_{1} m_{s}^{a}, && Q^{ab} = \kappa_{\star}^{2}q_{w}^{ab} + \kappa_{1}^2 q_{s}^{ab}.
\end{align}
\noindent where $\kappa_{1}=\mathbb{E}_{z}\left[z\sigma(z)\right]$ and $\kappa_{\star}^2=\mathbb{E}_{z}\left[\sigma(z)^2\right]-\kappa_{1}^2$ as in eq.~\eqref{eq:app:kappas}. With this notation, the asymptotic joint probability is simply written as:
\begin{align}
P(\nu_{\mu}, \{\lambda_{\mu}^{a}\}_{a=1}^{r}) = \frac{1}{\sqrt{\det\left(2\pi \Sigma\right)}} e^{-\frac{1}{2}\sum\limits_{a,b=0}^{r}z_{\mu}^{a}\left(\Sigma^{-1}\right)^{ab}z^{b}_{\mu}}
\end{align}
\noindent with $z_{\mu}^{0} = \nu_{\mu}$ and $z_{\mu}^{a} = \lambda^{a}_{\mu}$ for $a=1,\cdots, r$. The average over the replicated partition function (\ref{eq:Zr-av}) therefore reads:
\begin{align}
	\mathbb{E}_{\lbrace\vec{x}^{\mu},y^{\mu}\rbrace}\mathcal{Z}^{r}_{\beta} &= \int\dd\vec{\theta}^{0}~P_{\theta}(\vec{\theta}^{0})\int\prod\limits_{a=1}^{r}\dd\vec{w}~P_{w}\left(\vec{w}^{a}\right)\prod\limits_{\mu=1}^{n}\int\dd y^{\mu}\times\notag\\
	&\qquad\times\int\dd\nu_{\mu}~P_{y}^{0}\left(y^{\mu}|\nu_{\mu}\right)\int\prod\limits_{a=1}^{r}\dd \lambda^{a}_{\mu}~P(\nu_{\mu},\{\lambda_{\mu}^{a}\}) \prod\limits_{a=1}^{r}P_{y}\left(y^{\mu}|\{\lambda_{\mu}^a\}\right).
	\label{eq:app:avg_Z}
\end{align}

\subsection*{Rewriting as a saddle-point problem}
Note that after taking the average over $\vec{x}$, the integrals involved in the replicated partition function only couple through the overlap parameters. It is therefore useful to introduce the following Dirac $\delta$-functions to unconstrain them, introducing the decomposition:
\begin{align}
1 &= d^{-(r+1)^2} \int \dd \rho~\delta\left(d \rho-||\vec{\theta}^{0}||^2\right)\int \prod\limits_{a=1}^{r}\dd m_{s}^{a}~\delta\left(d m^{a}_{s}-\vec{s}^{a}\cdot \vec{\theta^{0}}\right)\times\notag\\
&\qquad\times\int\prod\limits_{1\leq a\leq b\leq r}\dd q_{s}^{ab}\delta\left(d q_{s}^{ab}-\vec{s}^{a}\cdot\vec{s}^{b}\right)\int\prod\limits_{1\leq a\leq b\leq r}\dd q_{w}^{ab}~\delta\left(p q_{w}^{ab}-\vec{w}^{a}\cdot\vec{w}^{b}\right)\notag\\
&= d^{-(r+1)^2} \int \frac{\dd \rho\dd\hat{\rho}}{2\pi}~e^{-i\hat{\rho}\left(d \rho-||\vec{\theta}^{0}||^2\right)}\int \prod\limits_{a=1}^{r}\frac{\dd m_{s}^{a}\dd \hat{m}_{s}^{a}}{2\pi}~e^{-i\sum\limits_{a=1}^{r}\hat{m}^{a}_{s}\left(d m^{a}_{s}-\vec{s}^{a}\cdot \vec{\theta^{0}}\right)}\times\notag\\
&\qquad\times\int\prod\limits_{1\leq a\leq b\leq r}\frac{\dd q_{s}^{ab}\dd \hat{q}_{s}^{ab}}{2\pi}e^{-i\sum\limits_{1\leq a\leq b\leq r}\hat{q}_{s}^{ab}\left(d q_{s}^{ab}-\vec{s}^{a}\cdot\vec{s}^{b}\right)}\int\prod\limits_{1\leq a\leq b\leq r}\frac{\dd q_{w}^{ab}\hat{q}^{ab}_{w}}{2\pi}~e^{-i\sum\limits_{1\leq a\leq b\leq r}\hat{q}_{w}^{ab}\left(p q_{w}^{ab}-\vec{w}^{a}\cdot\vec{w}^{b}\right)}.
\end{align}
Introducing the above in eq.~\eqref{eq:app:avg_Z} and exchanging the integration order allows to factorise the integrals over the $d,p,n$ dimensions and rewrite:
\begin{align}
	\mathbb{E}_{\lbrace\vec{x}^{\mu},y^{\mu}\rbrace}\mathcal{Z}^{r}_{\beta} = \int\frac{\dd \rho\dd\hat{\rho}}{2\pi}\int \prod\limits_{a=1}^{r}\frac{\dd m_{s}^{a}\dd \hat{m}_{s}^{a}}{2\pi}\int\prod\limits_{1\leq a\leq b\leq r}\frac{\dd q_{s}^{ab}\dd \hat{q}_{s}^{ab}}{2\pi}\frac{\dd q_{w}^{ab}\dd\hat{q}^{ab}_{w}}{2\pi}e^{p\Phi^{(r)}}
	\label{eq:app:saddle}
\end{align}
\noindent where the integrals over the variables $m_s^a$, $q_s^{ab}$ and $q_w^{ab}$ run over $\mathbb{R}$, while those over
$\hat m_s^a$, $\hat q_s^{ab}$ and $\hat q_w^{ab}$ run over $i \mathbb{R}$. The function
$\Phi^{(r)}$, a function of all the overlap parameters, is given by:
\begin{align}
\Phi^{(r)} &= -\gamma \rho\hat{\rho} - \gamma \sum\limits_{a=1}^{r}m^{a}_{s}\hat{m}^{a}_{s}-\sum\limits_{1\leq a\leq b\leq r}\left(\gamma q_{s}^{ab}\hat{q}^{ab}_{s}+q_{w}\hat{q}_{w}\right)+\alpha\Psi_{y}^{(r)}\left(\rho,m^{a}_{s},q^{ab}_{s},q_{w}^{ab}\right)\notag\\
&\qquad+\Psi^{(r)}_{w}\left(\hat{\rho},\hat{m}^{a}_{s},\hat{q}^{ab}_{s},\hat{q}_{w}^{ab}\right)
\end{align}
\noindent where we recall that $\alpha = n/p$, $\gamma = d/p$, and we have introduced:
\begin{align}
\Psi_{y}^{(r)} &= \log\int\dd y\int\dd\nu~P^{0}_{y}\left(y|\nu\right)\int\prod\limits_{a=1}^{r}\left[\dd\lambda^{a}  P_{y}\left(y|\lambda^{a}\right)\right] P(\nu,\left\{\lambda^{a}\right\})\notag\\
\Psi_{w}^{(r)} &= \frac{1}{p}\log\int\dd \vec{\theta}^{0}P_{\theta}(\vec{\theta}^{0})e^{-\hat{\rho}||\vec{\theta}^{0}||
^2}\int\prod\limits_{a=1}^{r}\dd\vec{w}^{a}~P_{w}(\vec{w}^{a})e^{\sum\limits_{1\leq a\leq b\leq r}\left[\hat{q}_{w}^{ab}\vec{w}^{a}\cdot \vec{w}^{b}+\hat{q}_{s}^{ab}\vec{s}^{a}\cdot\vec{s}^{b}\right]-\sum\limits_{a=1}^{r}\hat{m}^{a}_{s}\vec{s}^{a}\cdot \vec{\theta}^{0}}
\end{align}
Note that $\vec{s}^{a} = \frac{1}{\sqrt{p}}\mat{F}\vec{w}^{a}$ is a function of $\vec{w}^{a}$, and must be kept under the $\vec{w}^{a}$ integral. In the thermodynamic limit where $p\to\infty$ with $n/p$ and $d/p$ fixed, the integral in eq.~\eqref{eq:app:saddle} concentrates around the values of the overlap parameters that extremize $\Phi^{(r)}$, and therefore
\begin{align}
f=-\lim\limits_{r\to 0^{+}}\frac{1}{r}\underset{\substack{\{\rho, \hat{\rho}, m_{s}^{a},\hat{m}_{s}^{a}\}\\     \{ q_{s}^{ab}, \hat{q}_{s}^{ab},q_{w}^{ab}, \hat{q}_{w}^{ab}\}}}{\extr} ~\Phi^{(r)}.
\label{eq:app:free_prelimit}
\end{align}
\noindent
\subsection*{Replica symmetric Ansatz}
In order to proceed with the $r\to 0^{+}$ limit, we restrict the extremization above to the following replica symmetric Ansatz:
\begin{align}
m_{s}^{a} = m_{s} && \hat{m}^{a} = \hat{m}_{s} && \text{ for } a=1,\cdots, r\notag\\
q_{s/w}^{aa} = r_{s/w} && \hat{q}_{s/w}^{aa} = -\frac{1}{2}\hat{r}_{s/w} && \text{ for } a=1,\cdots, r\notag\\
q_{s/w}^{ab} = q_{s/w} && \hat{q}_{s/w}^{ab} = \hat{q}_{s/w} && \text{ for } 1\leq a < b \leq r
\end{align}
Note that, in the particular case of a convex loss function with $\lambda>0$, the replica symmetric Ansatz is justified: the problem only admitting one solution, it \emph{a fortiori} coincides with the replica symmetric one. For non-convex losses, solutions that are not replica symmetric (also known as \emph{replica symmetry breaking}) are possible, and the energy landscape of the free energy needs to be carefully analysed. In the practical applications explored in this manuscript, we focus on convex losses with ridge regularisation, and therefore the replica symmetric assumption is fully justified.


Before proceeding with the limit in eq.~\eqref{eq:app:free_prelimit}, we need to verify that the above Ansatz is well defined - in other words, that we have not introduced a spurious order one term in $\Phi$ that would diverge. This means we need to check that $\lim\limits_{r\to 0^{+}}\Phi = 0$. 

First, with a bit of algebra one can check that, within our replica symmetric Ansatz:
\begin{align}
\lim\limits_{r\to 0^{+}}\Psi_{y}^{(r)} = 0	.
\end{align}
Therefore,
\begin{align}
\lim\limits_{r\to 0^{+}}\Phi ^{(r)}  = -\gamma\rho\hat{\rho} +\gamma\log\int_{\mathbb{R}}\dd\theta^{0}~P_{\theta}\left(\theta^{0}\right) e^{\hat{\rho}{\theta^{0}}^{2}}
\end{align}
\noindent where we have used the fact that $P_{\theta}$ is a factorised distribution to take the $p\to\infty$ limit. In order for this limit to be $0$, we need that $\hat{\rho}=0$, which also fixes $\rho$ to be a constant given by the second moment of $\theta^{0}$:
\begin{align}
	\rho = \mathbb{E}_{\theta^{0}}\left[{\theta^{0}}^{2}\right]
\end{align}
We now proceed with the limit in eq.~\eqref{eq:app:free_prelimit}. Let's look first at $\Psi_{y}$. The non-trivial limit comes from the fact that $\det\Sigma$ and $\Sigma^{-1}$ are non-trivial functions of $r$. It is not hard to see, however, that $\Sigma^{-1}$ itself has replica symmetric structure, with components given by:
\begin{align}
    \left(\Sigma^{-1}\right)^{00} &= \tilde{\rho} = \frac{R+(r-1)Q}{\rho(R+(r-1)Q)-r M^2}, && \left(\Sigma^{-1}\right)^{aa}=\tilde{R} = \frac{\rho R+(r-2)\rho~Q-(r-1)M^2}{(R-Q)(\rho~R+(r-1)\rho~Q-r~M^2)}\notag\\
    \left(\Sigma^{-1}\right)^{a0} &= \tilde{M} = \frac{M}{r~M^2-\rho~R-(r-1)\rho~Q}, && \left(\Sigma^{-1}\right)^{ab} = \tilde{Q} = \frac{M^2 - \rho~Q}{(R-Q)(\rho~R+(r-1)\rho~Q-r~M^2)}\label{eq:app:inverseq}
\end{align}
\noindent where $M$, $Q$ and $R$ are the rescaled overlap parameters in the replica symmetric Ansatz, that is:
\begin{align}
M = \kappa_{1} m_{s}, && Q = \kappa_{\star}^{2}q_{w} + \kappa_{1}^2 q_{s}, && R = \kappa_{\star}^{2}r_{w} + \kappa_{1}^2 r_{s}.
\end{align}
\noindent This allows us to write:
\begin{align}
\Psi_{y}^{(r)} &= \log\int\dd y\int\dd\nu~P^{0}_{y}\left(y|\nu\right)e^{-\frac{\tilde{\rho}}{2}\nu^{2}}\int\prod\limits_{a=1}^{r}\dd\lambda^{a}P_{y}\left(y|\lambda^{a}\right)e^{-\frac{\tilde{Q}}{2}\sum\limits_{a,b=1}^{n}\lambda^{a}\lambda^{b}-\frac{\tilde{R}-\tilde{Q}}{2}\sum\limits_{a=1}^{r}\left(\lambda^{a}\right)^2-\tilde{M}\nu\sum\limits_{a=1}^{n}\lambda^{a}}\notag\\
&\qquad-\frac{1}{2}\log\det\left(2\pi\Sigma\right).
\end{align}
In order to completely factor the integral above in the replica space, we use the \emph{Hubbard-Stratonovich transformation}:
\begin{align}
e^{-\frac{\tilde{Q}}{2}\sum\limits_{a,b=1}^{r}\lambda^{a}\lambda^{b}} = \mathbb{E}_{\xi}e^{\sqrt{-\tilde{Q}}\xi\sum\limits_{a=1}^{r}\lambda^{a}}	
\end{align}
\noindent for $\xi\sim\mathcal{N}(0,1)$, such that
\begin{align}
\Psi_{y}^{(r)} &= \mathbb{E}_{\xi} \log\int\dd y\int\dd\nu~P^{0}_{y}\left(y|\nu\right)e^{-\frac{\tilde{\rho}}{2}\nu^{2}}\left[\int\dd\lambda P_{y}\left(y|\lambda\right)e^{-\frac{\tilde{R}-\tilde{Q}}{2}\lambda^2+\left(\sqrt{-\tilde{Q}}\xi-\tilde{M}\nu\right)\lambda}\right]^{r}\notag\\
&\qquad-\frac{1}{2}\log\det\left(2\pi\Sigma\right).
\end{align}
Taking into account the $r$ dependence of the inverse elements and of the determinant, we can take the limit to get:
\begin{align}
\lim\limits_{r\to 0^{+}}\frac{1}{r}\Psi_{y}^{(r)}	&=\mathbb{E}_{\xi}\int_{\mathbb{R}}\dd y\int\frac{\dd\nu}{\sqrt{{2\pi \rho~}}} P^{0}_{y}\left(y|\nu\right)e^{-\frac{1}{2\rho~}\nu^2}\log{\int\frac{\dd\lambda}{\sqrt{2\pi}}~P_{y}\left(y|\lambda\right)}e^{-\frac{1}{2}\frac{\lambda^2}{R-Q}+\left(\frac{\sqrt{Q-M^2/\rho~}}{R-Q}\xi+\frac{M/\rho~}{R-Q}\nu\right)\lambda}\notag\\
&\qquad-\frac{1}{2}\log\left(R-Q\right)-\frac{1}{2}\frac{Q}{R-Q}
\end{align}
Finally, making a change of variables and defining:
\begin{align}
\mathcal{Z}^{\cdot/0}_{y}(y;\omega,V) = \int\frac{\dd x}{\sqrt{2\pi V}}e^{-\frac{1}{2V}(x-\omega)^2}P^{\cdot/0}_{y}\left(y|x\right)
\label{eq:app:auxiliary_Z}	
\end{align}
\noindent allows us to rewrite the limit of $\Psi_{y}$ - which abusing notation we still denote $\Psi_{y}$ - as:
\begin{align}
\Psi_{y} = 	\mathbb{E}_{\xi}\left[\int_{\mathbb{R}}\dd y~\mathcal{Z}_{y}^0\left(y; \frac{M}{\sqrt{Q}}\xi, \rho~-\frac{M^2}{Q}\right)\log\mathcal{Z}_{y}\left(y;\sqrt{Q}\xi, R-Q\right)\right].
\label{eq:app:psi_y}
\end{align}
One can follow a very similar approach for the limit of $\Psi_{w}$, although in this case the limit is much simpler, since there is no $r$ dependence on the hat variables. The limit can be written as:
\begin{align}
\Psi_{w} = \lim\limits_{p\to\infty}\frac{1}{p}\mathbb{E}_{\xi,\eta, \theta^{0}}\log\int_{\mathbb{R}^{d}}\dd\vec{s}~P_{s}(\vec{s};\eta)e^{-\frac{\hat{V}_{s}}{2}||\vec{s}||^2+\left(\sqrt{\hat{q}_{s}}\xi\vec{1}_{d}+\hat{m}_{s}\vec{\theta}^{0}\right)^{\top}\vec{s}}
\end{align}
\noindent for $\xi,\eta\sim\mathcal{N}(0,1)$, and we have defined:
\begin{align}
P_{s}(\vec{s};\eta) = \int_{\mathbb{R}^{p}}\dd\vec{w}~P_{w}(\vec{w})e^{-\frac{\hat{V}_{w}}{2}||\vec{w}||^2+\sqrt{\hat{q}_{w}}\eta\vec{1}_{p}^{\top}\vec{w}}\delta\left(\vec{s}-\frac{1}{\sqrt{p}}\mat{F}\vec{w}\right)
\end{align}

and we have defined the shorthands $\hat{V}_{w} = \hat{r}_{w}+\hat{q}_{w}$ and $\hat{V}_{s} = \hat{r}_{s}+\hat{q}_{s}$.
\subsection*{Summary of the replica symmetric free energy density}
Summarising the calculation above, the replica symmetric free energy density reads:
\begin{align}
f &= \extr ~\Big\{-\frac{\gamma}{2}r_{s}\hat{r}_{s}-\frac{\gamma}{2}q_{s}\hat{q}_{s}+\gamma m_{s}\hat{m}_{s} - \frac{1}{2}r_{w}\hat{r}_{w}-\frac{1}{2}q_{w}\hat{q}_{w}\notag\\
&\qquad\qquad-\alpha\Psi_{y}(R,Q,M)-\Psi_{w}\left(\hat{r}_{s},\hat{q}_{s},\hat{m}_{s},\hat{r}_{w},\hat{q}_{w}\right)\Big\}
\label{eq:app:replica_final}
\end{align}
\noindent with $\alpha = \frac{n}{p}$, $\gamma = \frac{d}{p}$, and:
\begin{align} 
Q = \kappa_{1}^2 q_{s}+\kappa_{\star}^2q_{w}, && R = \kappa_{1}^2 r_{s}+\kappa_{\star}^2r_{w} && M = \kappa_{1}m_{s}.
\end{align}
The so-called potentials $(\Psi_{y}, \Psi_{w})$ are given by:
\begin{align}
	\Psi_{w} &= \lim\limits_{p\to\infty}\frac{1}{p}\mathbb{E}_{\xi,\eta, \theta^{0}}\log\int_{\mathbb{R}^{d}}\dd\vec{s}P_{s}(\vec{s};\eta)e^{-\frac{\hat{V}_{s}}{2}||\vec{s}||^2+\left(\sqrt{\hat{q}_{s}}\xi\vec{1}_{d}+\hat{m}_{s}\vec{\theta}^{0}\right)^{\top}\vec{s}}\label{eq:app:Psiw}\\
\Psi_{y} &= 	\mathbb{E}_{\xi}\left[\int_{\mathbb{R}}\dd y~\mathcal{Z}_{y}^0\left(y; \frac{M}{\sqrt{Q}}\xi, \rho~-\frac{M^2}{Q}\right)\log\mathcal{Z}_{y}\left(y;\sqrt{Q}\xi, R-Q\right)\right].
\end{align}
\noindent where:
\begin{align}
P_{s}(\vec{s};\eta) &= \int_{\mathbb{R}^{p}}\dd\vec{w}~ P_{w}(\vec{w})e^{-\frac{\hat{V}_{w}}{2}||\vec{w}||^2+\sqrt{\hat{q}_{w}}\eta\vec{1}_{p}^{\top}\vec{w}}
\delta\left(\vec{s}-\frac{1}{\sqrt{p}}\mat{F}\vec{w}\right)\notag\\
\mathcal{Z}^{\cdot/0}_{y}(y;\omega,V) &= \int\frac{\dd x}{\sqrt{2\pi V}}e^{-\frac{1}{2V}(x-\omega)^2}P^{\cdot/0}_{y}\left(y|x\right)	
\end{align}

\subsection{Evaluating $\Psi_{w}$ for ridge regularisation and Gaussian prior}
Note that as long as the limit in $\Psi_{w}$ is well defined, the eq.~\eqref{eq:app:replica_final} holds for any $P_{\theta}$ and $P_{w}$. However, as discussed in Sec.~\ref{sec:app:data}, we are interested in $\vec{\theta}^{0}\sim\mathcal{N}(0,\mat{I}_{d})$ and ridge regularisation so that $P_{w} = \exp\left(-\frac{\beta\lambda}{2}||\vec{w}||^2\right)$. In this case, we simply have:
\begin{align}
P(\vec{s};\eta) = \frac{e^{\frac{p}{2}\frac{\eta^2\hat{q}_{w}}{\beta\lambda+\hat{V}_{w}}}}{(\beta\lambda+\hat{V}_{w})^{p/2}}\mathcal{N}(\vec{s}; \vec{\mu}, \Sigma)	
\end{align}
\noindent with:
\begin{align}
\vec{\mu} = \frac{\sqrt{\hat{q}_{w}}\eta}{\beta\lambda+\hat{V}_{w}}\frac{\mat{F}\vec{1}_{p}}{\sqrt{p}}\in\mathbb{R}^{d}, && \Sigma = \frac{1}{\beta\lambda+\hat{V}_{w}} \frac{\mat{F}\mat{F}^{\top}}{p}\in\mathbb{R}^{d\times d}.
\end{align}
\noindent Therefore the argument of the logarithm in $\Psi_{w}$ is just another Gaussian integral we can do explicitly:
\begin{align}
\mathbb{E}_{s}e^{-\frac{\hat{V}_{s}}{2}||\vec{s}||^2+\vec{b}^{\top}\vec{s}} = \frac{e^{\frac{p}{2}\frac{\eta^2\hat{q}_{w}}{\beta\lambda+\hat{V}_{w}}}}{\left( \beta\lambda+\hat{V}_{w}\right)^{p/2}}\frac{e^{-\frac{1}{2}\vec{\mu}^{\top}\Sigma^{-1}\vec{\mu}+\frac{1}{2\hat{V}_{s}}||\vec{b}+\Sigma^{-1}\vec{\mu}||^2}}{\sqrt{\det\left(\mat{I}_{d}+\hat{V}_{s}\Sigma\right)}}e^{-\frac{1}{2\hat{V}_{s}}\left(b+\Sigma^{-1}\mu\right)^{\top}\left(\mat{I}_{d}+\hat{V}_{s}\Sigma\right)^{-1}\left(b+\Sigma^{-1}\mu\right)}
\end{align}
\noindent where we have defined the shorthand $\vec{b}=\left(\sqrt{\hat{q}_{s}}\xi\vec{1}_{d}+\hat{m}_{s}\vec{\theta}^{0}\right)\in\mathbb{R}^{d}$. Inserting back in eq.~\eqref{eq:app:Psiw} and taking the log,
\begin{align}
\Psi_{w}&=\lim\limits_{p\to\infty}\mathbb{E}_{\theta^{0},\xi,\eta}\left[\frac{1}{2}\frac{\eta^2\hat{q}_{w}}{\beta\lambda+\hat{V}_{w}}-\frac{1}{2}\log\left(\beta\lambda+\hat{V}_{w}\right)-\frac{1}{2p}\tr\log\left(\mat{I}_{d}+\hat{V}_{s}\Sigma\right)-\frac{1}{2p}\vec{\mu}^{\top}\Sigma^{-1}\vec{\mu}\right.\notag\\
&\hspace{3cm}\left.+\frac{1}{2p\hat{V}_{s}}||\vec{b}+\Sigma^{-1}\vec{\mu}||^2-\frac{1}{2p\hat{V}_{s}}\left(b+\Sigma^{-1}\mu\right)^{\top}\left(\mat{I}_{d}+\hat{V}_{s}\Sigma\right)^{-1}\left(b+\Sigma^{-1}\mu\right)\right]	
\end{align}
The averages over $\eta,\xi,\vec{\theta}^{0}$ simplify this expression considerably:
\begin{align}
\mathbb{E}_{\eta}\left[\vec{\mu}^{\top}\Sigma^{-1}\vec{\mu}\right]=\frac{1}{p}\frac{\hat{q}_{w}}{(\beta\lambda+\hat{V}_{w})^2}\left(\mat{F}\vec{1}_{p}\right)^{\top}\Sigma^{-1}\left(\mat{F}\vec{1}_{p}\right) = d \frac{\hat{q}_{w}}{\beta\lambda+\hat{V}_{w}}\notag\\
\mathbb{E}_{\eta,\xi,\theta^{0}}||b+\Sigma^{-1}\vec{\mu}||^2 =d(\hat{m}_{s}^{2}+\hat{q}_{s})+\frac{1}{p}\hat{q}_{w}\tr\left(\mat{F}\mat{F}^{\top}\right)^{-1}\notag
\end{align}
\begin{align}
\mathbb{E}_{\eta,\xi,\theta^{0}}\left(b+\Sigma^{-1}\vec{\mu}\right)^{\top}\left(\mat{I}_{d}+\hat{V}_{s}\Sigma\right)^{-1}\left(b+\Sigma^{-1}\vec{\mu}\right) &= \frac{1}{p}\hat{q}_{w}\tr\left[\mat{F}\mat{F}^{\top}\left(\mat{I}_{d}+\hat{V}_{s}\Sigma\right)^{-1}\right]\notag\\
&\qquad+(\hat{m}_{s}^{2}+\hat{q}_{s})\tr\left(\mat{I}_{d}+\hat{V}_{s}\Sigma\right)^{-1}
\end{align}
Finally, we can combine the two terms:
\begin{align}
\tr\frac{\mat{F}\mat{F}^{\top}}{p}+\tr\left[\frac{\mat{F}\mat{F}^{\top}}{p}\left(\mat{I}_{d}+\hat{V}_{s}\Sigma\right)^{-1}\right] = \frac{\hat{V}_{s}}{\beta\lambda+\hat{V}_{w}}\tr\left(\mat{I}_{d}+\hat{V}_{s}\Sigma\right)^{-1},
\end{align}
\noindent and write:
\begin{align}
\Psi_{w} &=-\frac{1}{2}\log\left(\beta\lambda+\hat{V}_{w}\right)-\frac{1}{2}\lim\limits_{p\to\infty}\frac{1}{p}\tr\log\left(\mat{I}_{d}+\frac{\hat{V}_{s}}{\beta\lambda+\hat{V}_{w}}\frac{\mat{F}\mat{F}^{\top}}{p}\right)\notag\\
&\qquad+\frac{\hat{m}_{s}^2+\hat{q}_{s}}{2\hat{V}_{s}}\left[\gamma-\lim\limits_{p\to\infty}\frac{1}{p}\tr\left(\mat{I}_{d}+\frac{\hat{V}_{s}}{\beta\lambda+\hat{V}_{w}}\frac{\mat{F}\mat{F}^{\top}}{p}\right)^{-1}\right]\notag\\
&\qquad+\frac{1}{2}\frac{\hat{q}_{w}}{\beta\lambda+\hat{V}_{w}}\left[1-\gamma+\lim\limits_{p\to\infty}\frac{1}{p}\tr\left(\mat{I}_{d}+\frac{\hat{V}_{s}}{\beta\lambda+\hat{V}_{w}}\frac{\mat{F}\mat{F}^{\top}}{p}\right)^{-1}\right]
\end{align}
Note that $\Psi$ only depends on the spectral properties of the matrix $\frac{1}{p}\mat{F}\mat{F}^{\top}\in\mathbb{R}^{d\times d}$, and more specifically on its resolvent in the asymptotic limit. A case of particular interest is when $\mat{F}\mat{F}^{\top}$ has a well defined spectral measure $\mu$ on the $p,d\to\infty$ limit with $\gamma = d/p$ fixed. In that case, we can write:
\begin{align}
	\lim\limits_{p\to\infty}\frac{1}{p}\tr\left(\mat{I}_{d}+\frac{\hat{V}_{s}}{\beta\lambda+\hat{V}_{w}}\frac{\mat{F}\mat{F}^{\top}}{p}\right)^{-1} &=\gamma\frac{\beta\lambda+\hat{V}_{w}}{\hat{V}_{s}}g_{\mu}\left(-\frac{\beta\lambda+\hat{V}_{w}}{\hat{V}_{s}}\right)\\
\end{align}
\noindent where $g_{\mu}$ is the Stieltjes transform of $\mu$, defined by:
\begin{align}
g_{\mu}(z) = \int\frac{\dd\mu(t)}{t-z}.	
\end{align}
Similarly, the logarithm term can be expressed as the logarithm potential of $\mu$ - although for the purpose of evaluating the generalisation error we will only need the derivative of these terms, and therefore only the Stieltjes transforms and its derivative.

In what follows, we will mostly focus on two kinds of projection matrices $\mat{F}$:
\paragraph{Gaussian projections: } For $\mat{F} \in\mathbb{R}^{d\times p}$ a random matrix with i.i.d. Gaussian entries with zero mean and variance $1$, $\mu$ is given by the well-known Marchenko-Pastur law, and the corresponding Stieltjes transform is given by:
\begin{align}
	g_{\mu}(z) = \frac{1-z-\gamma-\sqrt{(z-1-\gamma)^2-4\gamma}}{2z\gamma}, && z<0
\end{align}

\paragraph{Orthogonally invariant projection:} For $\mat{F} = \mat{U}^{\top}\mat{D}\mat{V}$ with $\mat{U}\in\mathbb{R}^{d\times d}$ and $\mat{V}\in\mathbb{R}^{p\times p}$ two orthogonal matrices and $\mat{D}\in\mathbb{R}^{d\times p}$ a rectangular diagonal matrix of rank $\min(d,p)$ and diagonal entries $d_{k}$, the empirical spectral density $\mu_{p}$ is given by:
\begin{align}
\mu_{d}(\lambda) = \frac{1}{d}\sum\limits_{k=1}^{\min(r,p)}\delta(\lambda-\lambda_{k})=\left(1-\min\left(1,\frac{1}{\gamma}\right)\right)\delta(\lambda)+\frac{1}{p}\sum\limits_{k=1}^{\min(d,p)}\delta(\lambda-d_{k}^2)
\end{align}
Therefore the choice of diagonal elements $d_{k}$ fully characterise the spectrum of $\mat{F}\mat{F}^{\top}$. In order for the orthogonally invariant case to be comparable to the Gaussian case, we fix $d_{k}$ in such a way that the projected vector $\mat{F}\vec{w}$ is of the same order in both cases, i.e.
\begin{align}
d_{k}^2	=
\begin{cases}
	\gamma & \text{ for }\gamma >1\\
	1 & \text{ for } \gamma \leq 1
\end{cases}
\end{align}
With this choice, the Stieltjes transform of $\mu$ reads:
\begin{align}
g_{\mu}(z) =
\begin{cases}
-(1-\frac{1}{\gamma})\frac{1}{z}+\frac{1}{\gamma}\frac{1}{\gamma-z} & \text{ for }\gamma >1\\
\frac{1}{1-z} & \text{ for }\gamma \leq 1\\
\end{cases}
\end{align}

\subsection{Gaussian equivalent model}
\label{sec:app:equiv_model}
It is interesting to note that the average over the dataset $\{\vec{x}^{\mu}, y^{\mu}\}_{\mu=1}^{n}$ of the replicated partition function $\mathcal{Z}_{\beta}^{r}$ in eq.~\eqref{eq:app:avg_Z}, obtained after the application of the GET, is identical to the replicated partition function of the same task over the following dual dataset $\{\vec{\tilde{x}}^{\mu},y^{\mu}\}_{\mu=1}^{n}$, where:
\begin{align}
\vec{\tilde{x}}^{\mu} = \kappa_{0}\vec{1}_{p}+\kappa_{1}	\frac{1}{\sqrt{d}}\mat{F}^{\top}\vec{c}^{\mu}+\kappa_{\star}\vec{z}^{\mu}
\end{align}
\noindent where $\vec{z}^{\mu}\sim\mathcal{N}(\vec{0}, \mat{I}_{p})$, and the labels $y^{\mu}\sim P_{y}$ are the same. Indeed, calling $\tilde{\mathcal{Z}}_{\beta}^{r}$ the replicated partition function for this equivalent dataset, and considering $\kappa_{0}$ we have:
\begin{align}
\mathbb{E}_{\{\vec{\tilde{x}}^{\mu},y^{\mu}\}}\tilde{\mathcal{Z}}_{\beta}^{r} &= \int\dd\vec{\theta}^{0}~P_{\theta}(\vec{\theta}^{0})\int\prod\limits_{a=1}^{r}\dd\vec{w}~P_{w}\left(\vec{w}^{a}\right)\times\notag\\
                                          &\qquad\times\prod\limits_{\mu=1}^{n}\int\dd y^{\mu}~\underbrace{\mathbb{E}_{\vec{c}^{\mu},\vec{z}^{\mu}}\left[P_{y}^{0}\left(y^{\mu}\big|\frac{\vec{c}^{\mu}\cdot\vec{\theta}^{0}}{\sqrt{d}}\right)\prod\limits_{a=1}^{r}P_{y}\left(y^{\mu}|\vec{w}^{a}\cdot \left(\frac{\kappa_{1}}{\sqrt{d}}\mat{F}^{\top}\vec{c}^{\mu}+\kappa_{\star}\vec{z}^{\mu}\right)  \right)\right]}_{\text{(I)}}.
                             \label{eq:Zr-avdual}
\end{align}
Rewriting (I):
\begin{align}
	\text{(I)} &= \int\dd\nu_{\mu} ~P_{y}^{0}\left(y^{\mu}|\nu_{\mu}\right)\int\prod\limits_{a=1}^{r}\dd\lambda_{\mu}^{a} ~P_{y}\left(y^{\mu}|\lambda^{a}_{\mu}\right)\times\notag\\
	&\qquad\times\underbrace{\mathbb{E}_{\vec{c}^{\mu},\vec{z}^{\mu}}\left[\delta\left(\nu_{\mu}-\frac{1}{\sqrt{d}}\vec{c}^{\mu}\cdot \vec{\theta}^{0}\right)\prod\limits_{a=1}^{r}\delta\left(\lambda_{\mu}^{a}-\frac{\kappa_{1}}{\sqrt{d}}\vec{w}^{a}\cdot \mat{F}^{\top}\vec{c}^{\mu}+\kappa_{\star}\vec{w}^{a}\cdot \vec{z}^{\mu}\right)\right]}_{\equiv P(\nu,\lambda)}.
\end{align}
It is easy to show that taking $(\kappa_{0},\kappa_{1})$ to match those from eq.~\eqref{eq:app:kappas}, the variables $\left(\nu_{\mu}, \{\lambda^{a}_{\mu}\}\right)$ are jointly Gaussian variables with correlation matrix given by $\Sigma$ exactly as in eq.~\eqref{eq:app:correlation}. This establishes the equivalence
\begin{align}
\tilde{\mathcal{Z}}_{\beta}^{r}	= \mathcal{Z}_{\beta}^{r}
\end{align}
\noindent from which follows the equivalence between the asymptotic generalisation and test error of these two models.

\section{Saddle-point equations and the generalisation error}
\label{sec:app:sp}
The upshot of the replica analysis is to exchange the $p$-dimensional minimisation problem for $\vec{w}\in\mathbb{R}^{p}$ in eq.~\eqref{eq:app:minimisation} for a one-dimensional minimisation problem for the parameters $\{r_{s}, q_{s}, m_{s}, r_{w}, q_{w}\}$ and their conjugate in eq.~\eqref{eq:app:replica_final}. In particular, note that by construction at the limit $\beta\to\infty$ the solution $\{q^{\star}_{s}, m^{\star}_{s}, q^{\star}_{w}\}$ of eq.~\eqref{eq:app:replica_final} corresponds to:
\begin{align}
	q_{w}^{\star} = \frac{1}{p}||\hat{\vec{w}}||^2 && q_{s}^{\star} = \frac{1}{d}||\mat{F}\hat{\vec{w}}||^2 && m_{s}^{\star} = \frac{1}{d}\left(\mat{F}\hat{\vec{w}}\right)\cdot \vec{\theta}^0
	\label{eq:app:optimal_overlaps}
\end{align}
\noindent where $\vec{\hat{w}}$ is the solution of the solution of eq.~\eqref{eq:app:minimisation}. As we will see, both the generalisation error defined in eq.~\eqref{eq:app:generalisation_error} and the training loss can be expressed entirely in terms of these overlap parameters.

\subsection{Generalisation error as a function of the overlaps} 
Let $\{\vec{x}^{\new}, y^{\new}\}$ be a new sample independently drawn from the same distribution of our data $\{\vec{x}^{\mu},y^{\mu}\}_{\mu=1}^{n}$. The generalisation error can then be written as:
\begin{align}
\epsilon_{g} &= \frac{1}{4^{k}}\mathbb{E}_{\vec{x}^{\new}, y^{\new}}\left(y^{\new}-\hat{f}\left(\sigma\left(\mat{F}^{\top}\vec{c}^{\new}\right)\cdot\vec{\hat{w}}\right)\right)^2\notag\\
&= \frac{1}{4^{k}}\int\dd y\int\dd \nu~P_{y}^0(y|\nu)\int\dd\lambda~(y-\hat{f}(\lambda))^2\mathbb{E}_{\vec{c}^{\new}}\left[\delta\left(\nu - \vec{c}^{\new}\cdot \vec{\theta}^0\right)\delta\left(\lambda - \sigma\left(\mat{F}^{\top}\vec{c}^{\new}\right)\cdot \vec{\hat{w}}\right)\right].
\label{eq:app:gen_error_overlaps}
\end{align}
\noindent where for convenience, we normalise $k=0$ for the regression task and $k=1$ for the classification task. Again, we apply the GET from Sec.~\ref{sec:app:get} to write the joint distribution over $\{\nu,\lambda\}$:
\begin{align}
P(\nu,\lambda) = \frac{1}{\sqrt{\det\left(2\pi \Sigma\right)}} e^{-\frac{1}{2}\vec{z}^{\top}\Sigma^{-1}\vec{z}},
\end{align}	
\noindent where $\vec{z} = (\nu,\lambda)^{\top}\in\mathbb{R}^{2}$ and $\Sigma$ is given by
\begin{align}
	\Sigma = 
	\begin{pmatrix}
		\rho & M^{\star} \\
		M^{\star} & Q^{\star}
	\end{pmatrix}, && \rho = \frac{1}{d}||\vec{\theta^0}||^2 && M^{\star} = \frac{\kappa_{1}}{d}\left(\mat{F}\vec{\hat{w}}\right)\cdot \vec{\theta}^0,  && Q^{\star} = \frac{\kappa_{1}^2}{d}||\mat{F}\vec{\hat{w}}||^2+\frac{\kappa_{\star}^{2}}{p}||\vec{\hat{w}}||^2.
\end{align}
Inserting in eq.~\eqref{eq:app:gen_error_overlaps} gives the desired representation of the generalisation error in terms of the optimal overlap parameters:
\begin{align}
	\epsilon_{g} =\frac{1}{4^{k}}\int\dd y\int\dd \nu~P_{y}^0(y|\nu)\int\dd\lambda~P(\nu,\lambda)(y-\hat{f}(\lambda))^2
\end{align}
For linear labels $y=\vec{c}\cdot\vec{\theta}^0$ in the regression problem, we simply have:
\begin{align}
\epsilon_{g} = \rho+Q^{\star}-2M^{\star}	
\end{align}
\noindent while for the corresponding classification problem with $y=\sign\left(\vec{c}\cdot\vec{\theta}^0\right)$:
\begin{align}
\epsilon_{g} = \frac{1}{\pi}\cos^{-1}\left(\frac{M^{\star}}{\sqrt{Q^{\star}}}\right)
\end{align}
\noindent which, as expected, only depend on the angle between $\mat{F}\vec{\hat{w}}$ and $\vec{\theta}^0$.

\subsection{Training loss}
Similarly to the generalisation error, the asymptotic of the training loss, defined for the training data $\{\vec{x}^{\mu},y^{\mu}\}_{\mu=1}^{n}$ as:
\begin{equation}
\epsilon_{t} = \frac{1}{n} \mathbb{E}_{\lbrace\vec{x}^{\mu},y^{\mu}\rbrace}\left[ \sum_{\mu = 1}^{n} \ell\left( y^{\mu}, \vec{x}^{\mu} \cdot \vec{\hat{w}}\right) + \frac{\lambda}{2} \Vert \vec{\hat{w}} \Vert_2^2 \right], 
\end{equation}
\noindent can also be written only in terms of the overlap parameters. Indeed, it is closely related to the free energy density defined in eq. ~\eqref{eq:app:def_freeen}. A close inspection on this definition tells us that:
\begin{equation}
\lim\limits_{n\to\infty}\epsilon_{t} = \lim\limits_{\beta\to\infty}\partial_{\beta} f_{\beta}.
\end{equation}
Taking the derivative of the free energy with respect to the parameter $\beta$ and recalling that $p = \alpha n$, we can then get:
\begin{equation}
\lim\limits_{n\to\infty}\epsilon_{t} = \frac{\lambda}{2\alpha} \lim\limits_{p\to\infty} \mathbb{E}_{\lbrace\vec{x}^{\mu},y^{\mu}\rbrace} \left[ \frac{\Vert \vec{\hat{w}} \Vert_2^2}{p} \right] - \lim\limits_{\beta\to\infty}\partial_{\beta} \Psi_y. 
\end{equation}
\noindent For what concerns the contribution of the regulariser, we simply note that in the limit of $p\to\infty$, the average concentrates around the overlap parameter $q_w^{\star}$. Instead, for what concerns the contribution of the loss function, we can start by explicitly taking the derivative with respect to $\beta$ of $\Psi_y$ in eq.~\eqref{eq:app:psi_y}, i.e.:
\begin{equation}
\partial_{\beta} \Psi_y = - \mathbb{E}_{\xi} \left[\int_{\mathbb{R}}\dd y~ \frac{\mathcal{Z}_{y}^0\left(y, \omega_0^{\star} \right)}{\mathcal{Z}_y \left(y, \omega_1^{\star}\right)} \int\frac{\dd x}{\sqrt{2\pi V_1^{\star}}}e^{-\frac{1}{2V_1^{\star}}(x-\omega_1^{\star})^2 - \beta \ell\left(y,x \right) } \ell\left(y,x \right)\right] ,
\end{equation}
\noindent with $\mathcal{Z}_y^{\cdot/0}$ defined in eq.~ \eqref{eq:app:auxiliary_Z}. At this point, as explained more in details in section \ref{sec:beta_to_infinity}, we can notice that in the limit of $\beta\to\infty$, it holds:
\begin{equation}
\lim\limits_{\beta\to\infty}\partial_{\beta} \Psi_y = - \mathbb{E}_{\xi}  \left[\int_{\mathbb{R}}\dd y~ \mathcal{Z}_{y}^0\left(y, \omega_0^{\star} \right) \ell\left(y,\eta \left(y, \omega_1^{\star} \right)  \right) \right],
\end{equation} 
\noindent with $\eta\left(y, \omega_1^{\star} \right)$ given in eq.~ \eqref{eq:app:eta}. 
Combining the two results together we then finally get:
\begin{equation}
\lim\limits_{n\to\infty} \epsilon_t \rightarrow \frac{\lambda}{2\alpha} q_w^{\star} + \mathbb{E}_{\xi} \left[\int_{\mathbb{R}}\dd y~  \mathcal{Z}_{y}^0\left(y, \omega_0^{\star} \right) \ell\left(y,\eta \left(y, \omega_1^{\star} \right) \right) \right]. 
\end{equation}

\subsection{Solving for the overlaps} 
As we showed above, both the generalisation error and the training loss are completely determined by the $\beta\to\infty$ solution of the extremization problem in eq.~\eqref{eq:app:replica_final}. For strictly convex losses $\ell$, there is a unique solution to this problem, that can be found by considering the derivatives of the replica potential. This leads to a set of self-consistent saddle-point equations that can be solved iteratively:
\begin{align}
\begin{cases}
	\hat{r}_{s} = -2\frac{\alpha \kappa_1^2}{\gamma}\partial_{r_{s}}\Psi_{y}\left(R,Q,M\right)\\
	\hat{q}_{s} = -2\frac{\alpha \kappa_1^2}{\gamma}\partial_{q_{s}}\Psi_{y}\left(R,Q,M\right)\\
	\hat{m}_{s} = \frac{\alpha \kappa_1}{\gamma}\partial_{m_{s}}\Psi_{y}\left(R,Q,M\right)\\\\
	\hat{r}_{w} = -2\alpha \kappa_{\star}^2\partial_{r_{w}}\Psi_{y}\left(R,Q,M\right)\\
	\hat{q}_{w} = -2\alpha \kappa_{\star}^2 \partial_{q_{w}}\Psi_{y}\left(R,Q,M\right)\\
\end{cases}&&
\begin{cases}
	r_{s} = -\frac{2}{\gamma}\partial_{\hat{r}_{s}}\Psi_{w}\left(\hat{r}_{s},\hat{q}_{s},\hat{m}_{s},\hat{r}_{w},\hat{q}_{w}\right)\\
	q_{s} = -\frac{2}{\gamma}\partial_{\hat{q}_{s}}\Psi_{w}\left(\hat{r}_{s},\hat{q}_{s},\hat{m}_{s},\hat{r}_{w},\hat{q}_{w}\right)\\
	m_{s} = \frac{1}{\gamma}\partial_{\hat{m}_{s}}\Psi_{w}\left(\hat{r}_{s},\hat{q}_{s},\hat{m}_{s},\hat{r}_{w},\hat{q}_{w}\right)\\\\
	r_{w} = -2\partial_{\hat{r}_{w}}\Psi_{w}\left(\hat{r}_{s},\hat{q}_{s},\hat{m}_{s},\hat{r}_{w},\hat{q}_{w}\right)\\
	q_{w} = -2\partial_{\hat{q}_{w}}\Psi_{w}\left(\hat{r}_{s},\hat{q}_{s},\hat{m}_{s},\hat{r}_{w},\hat{q}_{w}\right)\\
\end{cases}
\end{align}
In the case of a $\mat{F}$ with well-defined spectral density $\mu$, we can be more explicit and write:
\begin{align}
\begin{cases}
	V_{s} = \frac{1}{\hat{V}_{s}}\left(1-z~g_{\mu}(-z)\right) \\
	q_{s} = \frac{\hat{m}_{s}^2+\hat{q}_{s}}{\hat{V}_{s}^2}\left[1-2z g_{\mu}(-z)+z^2 g_{\mu}'(-z)\right]-\frac{\hat{q}_{w}}{(\beta\lambda+\hat{V}_{w})\hat{V}_{s}}\left[-z g_{\mu}(-z)+z^2 g'_{\mu}(-z)\right] \\
	m_{s} = \frac{\hat{m}_{s}}{\hat{V}_{s}}\left(1-z~g_{\mu}(-z)\right)	\\ \\
	V_{w} =\frac{\gamma}{\beta\lambda+\hat{V}_{w}}\left[\frac{1}{\gamma}-1+z g_{\mu}(-z)\right] \\
	q_{w} = \gamma\frac{\hat{q}_{w}}{(\beta\lambda+\hat{V}_{w})^2}\left[\frac{1}{\gamma}-1+z^2 g'_{\mu}(-z)\right]-\gamma\frac{\hat{m}_{s}^2+\hat{q}_{s}}{(\beta\lambda+\hat{V}_{w})\hat{V}_{s}}\left[-z g_{\mu}(-z)+z^2g'_{\mu}(-z)\right]
\end{cases}
\end{align}
\noindent where:
\begin{align}
V_{s/w} = r_{s/w} - q_{r/w} && \hat{V}_{s/w} = \hat{r}_{s/w} + \hat{q}_{r/w} && z = \frac{\beta\lambda+\hat{V}_{w}}{\hat{V}_{s}}	
\end{align}
We can also simplify slightly the derivatives of $\Psi_{y}$ without loosing generality by applying Stein's lemma, yielding:
\begin{align}
	\begin{cases}
		\hat{V}_{s} = - \frac{\alpha\kappa_{1}^{2}}{\gamma} \mathbb{E}_{\xi}\left[\int_{\mathbb{R}}\dd y~\mathcal{Z}^{0}_{y}\left(y;\frac{M}{\sqrt{Q}}\xi, \rho-\frac{M^2}{Q}\right)\partial_{\omega}f_{y}\left(y;\sqrt{Q}\xi, R-Q\right)\right]\\
		\hat{q}_{s} = \frac{\alpha\kappa_{1}^2}{\gamma} \mathbb{E}_{\xi}\left[\int_{\mathbb{R}}\dd y~\mathcal{Z}^{0}_{y}\left(y;\frac{M}{\sqrt{Q}}\xi, \rho-\frac{M^2}{Q}\right)f_{y}\left(y;\sqrt{Q}\xi, R-Q\right)^2\right]\\
		\hat{m}_{s} = \frac{\alpha\kappa_{1}}{\gamma} \mathbb{E}_{\xi}\left[\int_{\mathbb{R}}\dd y~\mathcal{Z}^{0}_{y}\left(y;\frac{M}{\sqrt{Q}}\xi, \rho-\frac{M^2}{Q}\right)f^{0}_{y}\left(y;\frac{M}{\sqrt{Q}}\xi, \rho-\frac{M^2}{Q}\right)f_{y}\left(y;\sqrt{Q}\xi, R-Q\right)\right]\\ \\
		\hat{V}_{w} = - \alpha\kappa_{\star}^2 \mathbb{E}_{\xi}\left[\int_{\mathbb{R}}\dd y~\mathcal{Z}^{0}_{y}\left(y;\frac{M}{\sqrt{Q}}\xi, \rho-\frac{M^2}{Q}\right)\partial_{\omega}f_{y}\left(y;\sqrt{Q}\xi, R-Q\right)\right]\\
		\hat{q}_{w} = \alpha \kappa_{\star}^2\mathbb{E}_{\xi}\left[\int_{\mathbb{R}}\dd y~\mathcal{Z}^{0}_{y}\left(y;\frac{M}{\sqrt{Q}}\xi, \rho-\frac{M^2}{Q}\right)f_{y}\left(y;\sqrt{Q}\xi, R-Q\right)^2\right]
	\end{cases}
\end{align}
\noindent with $f^{\cdot/0}_{y}(y;\omega,V) = \partial_{\omega}\log\mathcal{Z}_{y}^{\cdot/0}$. For a given choice of spectral density $\mu$ (corresponding to a choice of projection $\mat{F}$), label rule $P^{0}_{y}$ and loss function $\ell$, the auxiliary functions $(\mathcal{Z}^{0},\mathcal{Z})$ can be computed, and from them the right-hand side of the update equations above. The equations can then be iterated until the convergence to the fixed point minimising the free energy at fixed $(\alpha,\gamma,\beta)$. For convex losses and $\beta\to\infty$, the fixed point of these equations gives the overlap corresponding to the estimator solving eq.~\eqref{eq:app:minimisation}.

\subsection{Taking $\beta\to\infty$ explicitly}
\label{sec:beta_to_infinity} 
Although the saddle-point equations above can be iterated explicitly for any $\beta>0$, it is envisageable to take the limit $\beta\to\infty$ explicitly, since $\beta$ is an auxiliary parameter we introduced, and that was not present in the original problem defined in eq.~\eqref{eq:app:minimisation}.

Since the overlap parameters depend on $\beta$ only implicitly through $\mathcal{Z}_{y}$ and its derivatives, we proceed with the following ansatz for their $\beta\to\infty$ scaling:
\begin{align}
V^{\infty}_{s/w}=\beta V_{s/w} && q_{s/w}^{\infty}= q_{s/w} && m_{s}^{\infty}= m_{s}\notag\\
\hat{V}^{\infty}_{s/w}=\frac{1}{\beta} \hat{V}_{s/w} && \hat{q}^{\infty}_{s/w} = \frac{1}{\beta^2}\hat{q}_{s/w} && \hat{m}^{\infty}_{s} = \hat{m}_{s}.
\label{eq:app:ansatz}
\end{align}
This ansatz can be motivated as follows. Recall that:
\begin{align}
	\mathcal{Z}_{y}(y;\omega,V) = \int\frac{\dd x}{\sqrt{2\pi V}}e^{-\beta\left[\frac{(x-\omega)^2}{2\beta V}+\ell(x,y)\right]} = \int\frac{\dd x}{\sqrt{2\pi V}}e^{-\beta\mathcal{L}(x)}.
	\label{eq:app:scalings}
\end{align}
Therefore, letting $V = \mu_{1}^{2}V_{s}+\mu_{\star}^{2}V_{w}$ scale as $V^{\infty}=\beta V$, at $\beta\to\infty$:
\begin{align}
	\mathcal{Z}_{y}(y;\omega,V) \underset{\beta\to\infty}{=} e^{-\beta \mathcal{L}(\eta)}
	\end{align}
\noindent where:
\begin{align}
\eta(y;\omega, V) = \underset{x\in\mathbb{R}}{\argmin}\left[\frac{(x-\omega)^2}{2 V^{\infty}}+\ell(x,y)\right].
\label{eq:app:eta}
\end{align}
For convex losses $\ell$ with $\lambda > 0$, this one-dimensional minimisation problem has a unique solution that can be easily evaluated. Intuitively, this ansatz translates the fact the variance of our estimator goes to zero as a power law at $\beta\to\infty$, meaning the Gibbs measure concentrates around the solution of the optimisation problem eq.~\eqref{eq:app:minimisation}. The other scalings in eq.~\eqref{eq:app:scalings} follow from analysing the dependence of the saddle-point equations in $V$. 

The ansatz in eq.~\eqref{eq:app:ansatz} allow us to take the $\beta\to\infty$ and rewrite the saddle-point equations as:
\begin{align}
	&\begin{cases}
		\hat{V}^{\infty}_{s} = \frac{\alpha\mu_{1}^{2}}{\gamma} \mathbb{E}_{\xi}\left[\int_{\mathbb{R}}\dd y~\mathcal{Z}^{0}_{y} \left(\frac{1-\partial_{\omega}\eta}{V^{\infty}}\right)\right]\\
		\hat{q}^{\infty}_{s} = \frac{\alpha\mu_{1}^2}{\gamma} \mathbb{E}_{\xi}\left[\int_{\mathbb{R}}\dd y~\mathcal{Z}^{0}_{y} \left( \frac{\eta	-\omega}{V^{\infty}}\right)^2 \right]\\
		\hat{m}^{\infty}_{s} = \frac{\alpha\mu_{1}}{\gamma} \mathbb{E}_{\xi}\left[\int_{\mathbb{R}}\dd y~\partial_{\omega}\mathcal{Z}^{0}_{y}\left(\frac{\eta-\omega}{V^{\infty}}\right) \right]\\ \\
		\hat{V}^{\infty}_{w} = \alpha\mu_{\star}^2 \mathbb{E}_{\xi}\left[\int_{\mathbb{R}}\dd y~\mathcal{Z}^{0}_{y}\left(\frac{1-\partial_{\omega}\eta}{V^{\infty}}\right)\right]\\
		\hat{q}^{\infty}_{w} = \alpha \mu_{\star}^2\mathbb{E}_{\xi}\left[\int_{\mathbb{R}}\dd y~\mathcal{Z}^{0}_{y}\left( \frac{\eta	-\omega}{V^{\infty}}\right)^2\right]
	\end{cases}\label{eq:app:eq_channel}
 \\ \notag\\
	&\begin{cases}
	V^{\infty}_{s} = \frac{1}{\hat{V}^{\infty}_{s}}\left(1-z~g_{\mu}(-z)\right) \\
	q^{\infty}_{s} = \frac{\left(\hat{m}^{\infty}_{s}\right)^2+\hat{q}^{\infty}_{s}}{\left(\hat{V}^{\infty}_{s}\right)^2}\left[1-2z g_{\mu}(-z)+z^2 g_{\mu}'(-z)\right]-\frac{\hat{q}^{\infty}_{w}}{(\lambda+\hat{V}_{w})\hat{V}_{s}}\left[-z g_{\mu}(-z)+z^2 g'_{\mu}(-z)\right] \\
	m^{\infty}_{s} = \frac{\hat{m}^{\infty}_{s}}{\hat{V}^{\infty}_{s}}\left(1-z~g_{\mu}(-z)\right)	\\ \\
	V^{\infty}_{w} =\frac{\gamma}{\lambda+\hat{V}^{\infty}_{w}}\left[\frac{1}{\gamma}-1+z g_{\mu}(-z)\right] \\
	q^{\infty}_{w} = \gamma\frac{\hat{q}^{\infty}_{w}}{(\lambda+\hat{V}^{\infty}_{w})^2}\left[\frac{1}{\gamma}-1+z^2 g'_{\mu}(-z)\right]-\gamma\frac{\left(\hat{m}^{\infty}_{s}\right)^2+\hat{q}^{\infty}_{s}}{(\lambda+\hat{V}^{\infty}_{w})\hat{V}^{\infty}_{s}}\left[-z g_{\mu}(-z)+z^2g'_{\mu}(-z)\right]
\end{cases}
\label{eq:app:eq_prior}
\end{align}
\noindent where $\mathcal{Z}_{y}^{0}(y;\omega,V)$ is always evaluated at $(\omega,V) = \left(\frac{M^{\infty}}{\sqrt{Q^{\infty}}}\xi, \rho-\frac{{M^{\infty}}^2}{Q^{\infty}}\right)$, $\eta(y;\omega, V)$ at $\left(\omega, V\right)=\left(\sqrt{Q^{\infty}}\xi, V^{\infty}\right)$ and $z = \frac{\lambda+\hat{V}^{\infty}_{w}}{\hat{V}^{\infty}_{s}}$.

\subsection{Examples}
In this section we exemplify our general result in two particular cases for which the integrals in the right-hand side of eq.~\eqref{eq:app:eq_channel} can be analytically performed: the ridge regression task with linear labels and a classification problem with square loss and ridge regularisation term. The former example appears in Fig.~\ref{fig:generalisation:erf:ridge:optimal} (left) and the later in Figs.~\ref{fig:generalisation:logistic:lambdas} (blue curve), \ref{fig:generalisation:logistic:HMM}, \ref{fig:generalisation:logistic:HMM_heat} of the main.

\paragraph{Ridge regression with linear labels:} Consider the task of doing ridge regression $\ell(y,x) =  \frac{1}{2}\left(y-x\right)^2$, $\lambda>0$ on the linear patterns $\vec{y} = \frac{1}{\sqrt{d}}\mat{C}\vec{\theta}^0+\sqrt{\Delta}\vec{z}$, with $\vec{z} \sim\mathcal{N}(\vec{0},\mat{I}_{n})$ and $\vec{\theta}^{\star}\sim\mathcal{N}(\vec{0},\mat{I}_{d})$. In this case, we have:
\begin{align}
\eta(y;\omega,V) = \frac{\omega+yV}{1+V}	\label{eq:app:xstar:square}
\end{align}
\noindent and the saddle-point equations for the hat overlap read:
\begin{align}
\hat{V}^{\infty}_{s} = \frac{\alpha}{\gamma}\frac{\kappa_{1}^2}{1+V^{\infty}}	&& \hat{q}^{0}_{s} = \frac{\alpha\kappa_{1}^2}{\gamma}\frac{1+\Delta+Q^{\infty}-2M^{\infty}}{\left(1+V^{\infty}\right)^2} && \hat{m}_{s} = \frac{\alpha}{\gamma}\frac{\kappa_{1}}{1+V^{\infty}}\notag
\end{align}
\begin{align}
	\hat{V}^{\infty}_{w} = \frac{\alpha\kappa_{\star}^2}{1+V^{\infty}}	&& \hat{q}^{\infty}_{w} = \alpha\kappa_{\star}^2\frac{1+\Delta+Q^{\infty}-2M^{\infty}}{\left(1+V^{\infty}\right)^2}
\end{align}
This particular example corresponds precisely to the setting studied in \cite{montanari2019generalisation}. 

\paragraph{Classification with square loss and ridge regularisation:} Consider a classification task with square loss $\ell(y,x) =  \frac{1}{2}\left(y- x\right)^2$ and labels generated as $\vec{y} = \sign\left(\frac{1}{\sqrt{d}}\mat{C}\vec{\theta}^{0}\right)$, with $\vec{\theta}^{0}\sim\mathcal{N}(\vec{0},\mat{I}_{d})$. Then the saddle-point equations are simply:
 \begin{align}
\hat{V}^{\infty}_{s} = \frac{\alpha}{\gamma}\frac{\kappa_{1}^2}{1+V^{\infty}}	&& \hat{q}^{\infty}_{s} = \frac{\alpha}{\gamma}\kappa_{1}^2\frac{1+Q^{\infty}-\frac{2\sqrt{2}M^{\infty}}{\sqrt{\pi}}}{\left(1+V^{\infty}\right)^2} && \hat{m}_{s} =\frac{\alpha}{\gamma}\sqrt{\frac{2}{\pi}}\frac{\kappa_{1}}{1+V^{\infty}}\notag
\end{align}
\begin{align}
	\hat{V}^{\infty}_{w} = \frac{\alpha\kappa_{\star}^2}{1+V^{\infty}}	&& \hat{q}^{\infty}_{w} = \alpha\kappa_{\star}^2\frac{1+Q^{\infty}-\frac{2M^{\infty}}{\sqrt{\pi}}}{\left(1+V^{\infty}\right)^2}
\end{align}

\section{Numerical Simulations}
\label{sec:app:sim}
In this section, we provide more details on how the numerical simulations in the main manuscript have been performed. 

First, the dataset $\{\vec{x}^{\mu}, y^{\mu}\}_{\mu=1}^{n}$ is generated according to the procedure described in Section \ref{sec:model} of the main, which we summarise here for convenience in algorithm \ref{alg:app:dataset}:

\begin{algorithm}[h]
   \caption{Generating dataset $\{\vec{x}^{\mu}, y^{\mu}\}_{\mu=1}^{n}$}
   \label{alg:app:dataset}
\begin{algorithmic}
   \STATE {\bfseries Input:} Integer $d$, parameters $\alpha,\gamma\in\mathbb{R}_{+}$, matrix $\mat{F}\in\mathbb{R}^{d\times p}$, vector $\vec{\theta}^{0}\in\mathbb{R}^{d}$ non-linear functions $\sigma, f^{0}:\mathbb{R}\to\mathbb{R}$.
   \STATE Assign $p \leftarrow \lfloor d/\gamma \rfloor$, $n \leftarrow \lfloor \alpha p \rfloor$
   \STATE Draw $\mat{C}\in\mathbb{R}^{n\times d}$ with entries $c_{\rho}^{\mu}\sim\mathcal{N}(0,1)$ i.i.d.
   \STATE Assign $\vec{y} \leftarrow f^{0}\left(\mat{C}\vec{\theta}^{0}\right)\in\mathbb{R}^{n}$ component-wise.
	\STATE Assign $\mat{X} \leftarrow \sigma\left(\mat{C}\mat{F}\right)\in\mathbb{R}^{n\times p}$ component-wise.
	\STATE {\bfseries Return:} $\mat{X}, \vec{y}$
\end{algorithmic}
\end{algorithm}
In all the examples from the main, we have drawn $\vec{\theta}^{0}\sim\mathcal{N}(0,\mat{I}_{d})$. For the regression task in Fig.~\ref{fig:generalisation:erf:ridge:optimal} we have taken $f^{0}(x) = x$, while in the remaining classification tasks $f^{0}(x) = \sign(x)$. For Gaussian projections, the components of $\mat{F}$ are drawn from $\mathcal{N}(0,1)$ i.i.d., and in for the random orthogonal projections we draw two orthogonal matrices $\mat{U}\in\mathbb{R}^{d\times d}$, $\mat{V}\in\mathbb{R}^{p\times p}$ from the Haar measure and we let $\mat{F} = \mat{U}^{\top}\mat{D}\mat{V}$ with $\mat{D}\in\mathbb{R}^{d\times p}$ a diagonal matrix with diagonal entries $d_{k} = \max(\sqrt{\gamma},1)$, $k=1,\cdots,\min(n,p)$.

Given this dataset, the aim is to infer the configuration $\hat{\vec{w}}$, minimising a given loss function with a ridge regularisation term. In the following, we describe how to accomplish this task for both square and logistic loss.

\paragraph{Square Loss:}
In this case, the goal is to solve the following optimisation problem:
\begin{equation}
\hat{\vec{w}} = \underset{\vec{w}}{\min}~ \left[\frac{1}{2}\sum\limits_{\mu=1}^{n} \left( y^{\mu} - \vec{x}^{\mu} \cdot \vec{w}\right)^2  +\frac{\lambda}{2}||\vec{w}||_{2}^2\right] \ .
\label{eq:app:ridge_sol}
\end{equation}
\noindent which has a simple closed-form solution given in terms of the Moore-Penrose inverse:
\begin{align}
\hat{\vec{w}} = \begin{cases}
 \left( \mat{X}^{\top}\mat{X} + \lambda \mat{I}_{p} \right)^{-1} \mat{X}^{\top} \vec{y} , & \text{if} \ n > p\\\\
 \mat{X}^{\top}\left(\mat{X}\mat{X}^T + \lambda \mat{I}_{n} \right)^{-1} \vec{y} , & \text{if} \ p > n 
 \end{cases}
\end{align}

\paragraph{Logistic Loss:}
In this case, the goal is to solve the following optimisation problem:
\begin{equation}
\hat{\vec{w}} = \underset{\vec{w}}{\min}~ \left[\sum\limits_{\mu=1}^{n} \log\left( 1 + e^{-y^{\mu}\left(\vec{x}^{\mu} \cdot \vec{w} \right)} \right)  +\frac{\lambda}{2}||\vec{w}||_{2}^2\right] \ .
\label{eq:app:logistic_sol}
\end{equation}
To solve the above, we use the \emph{Gradient Descent} (GD) on the regularised loss. In our simulations, we took advantage of Scikit-learn 0.22.1, an out-of-the-box open source library for machine learning tasks in Python \cite{scikit-learn, sklearn_api}. The library provides the class \emph{sklearn.linear\_{}model.LogisticRegression}, which implements GD with logistic loss and a further $\ell_2$-regularisation, if the parameter 'penalty' is set to 'l2'. GD stops either if the following condition is satisfied:
\begin{equation}
\max \{\left(  \nabla \vec{w}\right)_i  \vert i=1,...,p \} \leqslant \mbox{tol},
\end{equation} 
\noindent with $\nabla \vec{w}$ being the gradient, or if a maximum number of iterations is reached. We set tol to $10^{-4}$ and the maximum number of iterations to $10^{4}$.

In both cases described above, the algorithm returns the estimator $\hat{\vec{w}}\in\mathbb{R}^{p}$, from which all the quantities of interest can be evaluated. For instance, the generalisation error can be simply computed by drawing a new and independent sample $\{\mat{X}^{\text{new}}, \vec{y}^{\text{new}}\}$ using algorithm \ref{alg:app:dataset} with the same inputs $\mat{F}, \sigma, f^{0}$ and $\vec{\theta}^{0}$ and computing:
\begin{align}
\epsilon_{g}(n,p,d) = \frac{1}{4^{k} n}||\vec{y}^{\text{new}}-\hat{f}\left(\mat{X}^{\text{new}}\vec{\hat{w}}\right)||^{2}_{2}
\end{align}
\noindent with $\hat{f}(x)=x$ for the regression task and $\hat{f}(x) = \sign(x)$ for the classification task.

The procedure outlined above is repeated $n_{\text{seeds}}$ times, for different and independent draws of the random quantities $\mat{F}, \vec{\theta}^{0}$, and a simple mean is taken in order to obtain the ensemble average of the different quantities. In most of the examples from the main, we found that $n_{\text{seeds}}=30$ was enough to obtain a very good agreement with the analytical prediction from the replica analysis. The full pipeline for computing the averaged generalisation error is exemplified in algorithm \ref{alg:app:pipeline}.

\begin{algorithm}[h]
   \caption{Averaged generalisation error.}
   \label{alg:app:pipeline}
\begin{algorithmic}
   \STATE {\bfseries Input:} Integer $d$, parameters $\alpha,\gamma, \lambda\in\mathbb{R}_{+}$, non-linear functions $\sigma, f^{0}, \hat{f}$ and integer $n_{\text{seeds}}$.
   \STATE Assign $p \leftarrow \lfloor d/\gamma \rfloor$, $n \leftarrow \lfloor \alpha p \rfloor$
   \STATE Initialise $E_{g} = 0$.
   \FOR{$i=1$ {\bfseries to} $n_{\text{seeds}}$}
   \STATE Draw $\mat{F}$, $\vec{\theta}^{0}$.
   \STATE Assign $\mat{X}, \vec{y} \leftarrow$  Alg. \ref{alg:app:dataset}.
   \STATE Compute $\vec{\hat{w}}$ from eq.~\eqref{eq:app:ridge_sol} or \eqref{eq:app:logistic_sol} with $\mat{X},\vec{y}$ and $\lambda$.
   \STATE Generate new dataset $\mat{X}^{\text{new}}, \vec{y}^{\text{\new}}$ from Alg. \ref{alg:app:dataset}.
   \STATE Assign $E_{g} \leftarrow E_{g} + \frac{1}{4^{k} n}||\vec{y}^{\text{new}}-\hat{f}\left(\mat{X}^{\text{new}}\vec{\hat{w}}\right)||^{2}_{2} $
	\ENDFOR
	\STATE {\bfseries Return:} $\epsilon_{g} = \frac{E_{g}}{n_{\text{seeds}}}$
	\end{algorithmic}
\end{algorithm}

\newpage
\bibliographystyle{unsrt}
\bibliography{refs}

\end{document}